\documentclass[fleqn,twoside]{article}
\usepackage{espcrc1}

\newcommand{\AmS}{{\protect\the\textfont2
  A\kern-.1667em\lower.5ex\hbox{M}\kern-.125emS}}
\title{Composite Strategy for Multicriteria Ranking/Sorting
 (methodological issues, examples)}

\author{Mark Sh. Levin
\thanks{
 Mark Sh. Levin:~
 http://www.mslevin.iitp.ru;~
 email: mslevin@acm.org}
  }

\begin{document}

\maketitle

\begin{abstract}
 The paper addresses the modular design of
 composite solving strategies for multicriteria ranking
 (sorting).
 Here a ``scale of creativity''
 that is close to creative levels proposed by Altshuller
 is used as the reference viewpoint:
 (i) a basic object,
 (ii) a selected object,
 (iii) a modified object, and
 (iv) a designed object (e.g., composition of object components).
 These levels maybe used in various parts of
 decision support systems (DSS)
 (e.g., information, operations, user).
 The paper focuses on the more creative above-mentioned level
 (i.e., composition or combinatorial synthesis)
 for the operational part (i.e., composite solving strategy).
 This is important for a search/exploration mode of
 decision making process with usage
 of various procedures and techniques and
  analysis/integration
  of obtained results.

 The paper describes  methodological issues
  of decision technology
  and synthesis of composite strategy for multicriteria ranking.
 The synthesis of composite strategies is based on
 ``hierarchical morphological multicriteria design'' (HMMD)
 which
 is based on selection and combination of
 design alternatives (DAs)
 (here: local procedures or techniques)
 while taking into account
 their quality and quality of their
  interconnections (IC).
 A new version of HMMD with interval multiset estimates
 for DAs is used.
      The operational environment of DSS COMBI for
 multicriteria ranking, consisting of a morphology of
 local procedures or techniques
 (as design alternatives DAs),
 is examined as a basic one.

~~~~~~~~~~~

 {\it Keywords:}~
 Decision making technology;
  Multicriteria ranking; Sorting;
 Decision support system;
 System architecture,
  Problem solving strategies;
 Combinatorial synthesis;
  Software engineering;
   Multiset

\vspace{1pc}
\end{abstract}

\tableofcontents

\newcounter{cms}
\setlength{\unitlength}{1mm}

\section{INTRODUCTION}


 In recent years, there has been considerable interest in complex
 DSS
 consisting of the
  data bases (including descriptions of decision making situations),
  model bases,
  user interface, and
  knowledge bases
 (e.g., \cite{alavi92},\cite{bel92},\cite{bol02},\cite{enge03},\cite{eom93},\cite{eom95},\cite{eom98},
 \cite{for91},\cite{iij96},\cite{kau94},\cite{liu88},\cite{makov00},
 \cite{marak99},\cite{michal03},\cite{mora03},\cite{ribe03},\cite{rose97},\cite{saw87},
 \cite{sax93},\cite{seo93},\cite{sis99}, \cite{tur90}).
   The structures of multicriteria DSS were described in
 (\cite{ari85},\cite{ban93},\cite{bel92},\cite{da89},\cite{eom93},\cite{eom98},\cite{for03},
 \cite{how88},\cite{jel86},\cite{lev86},\cite{makov00},
 \cite{saw87},\cite{sis99},\cite{spr82},\cite{tur90}).
 Jelassi has proposed five generations of multicriteria DSS
 (\cite{jel86},\cite{sis99}).

 Main tendencies of decision support technology are presented in
  (\cite{enge03},\cite{sax93},\cite{shi02}).
    Saxena investigated computer-aided decision support engineering
 including six stages methodology for DSS development (problem definition,
 task analysis, requirements engineering, system design, system prototype,
 user evaluation and adaptation \cite{sax93}.
  Four models for a decision support system are
 analyzed in \cite{mir99}:
  (i) the symbolic DSS,
  (ii) the expert system,
  (iii) the holistic DSS, and
  (iv) the adaptive DSS.
   Adaptive decision support systems have been examined in
 (\cite{faz97},\cite{hols93}).
 Howard has described the decision analysis process \cite{how88}.
 An object relational approach for the design of DSS
  has been proposed by Srinivasan and Sundaram \cite{srin00}.
 A special framework to support the design of a DSS for some new
 decision problem was suggested in \cite{reich05}.
 Some contemporary trends are targeted to the following:
 (a) distributed cooperative DSS
 (e.g., \cite{bose96},\cite{gach03},\cite{shi02},\cite{zha02})
 (b) spatial DSS (SDDS, i.e., integration of DSS and GIS)
 (e.g.,
 \cite{cross95},\cite{parker03},\cite{rao07},\cite{rinner03});
 (c) multi-agent DSS (e.g.,
  \cite{bose96},\cite{bui99},\cite{hahn09},\cite{oss05});
 (d) mobile DSS
 (e.g.,
 \cite{eren08},\cite{heijden06},\cite{michal03},\cite{perez10});
 and
 (e) Web-based DSS
 (e.g.,
 \cite{bhar07},\cite{kwon03},\cite{rinner03},\cite{sal04},\cite{xie05},\cite{zhang07}).
 A classification of DSS is described in \cite{power04}.

 Issues of model management have been widely studied and used in DSS engineering
 (e.g.,
 \cite{app86},\cite{bin95},\cite{bla93},\cite{bol02},
 \cite{cha02},\cite{da89},\cite{dol84},\cite{gro91},\cite{hol97},\cite{jar88},\cite{lia85}).
 The research domain involves various directions, for example:
  ~{\it 1.} model selection (e.g., \cite{mang03}), e.g.,
  selection of a package for multi-attribute decision making
  that is more compatible with the user's needs \cite{french05};
 ~{\it 2.} usage of model libraries (e.g., \cite{man90});
 ~{\it 3.} building (construction, composition, integration)
 of models
 (e.g., \cite{bas94},\cite{basu94a},\cite{cha03},\cite{iij96},\cite{kri91});
 ~{\it 4.} an object-oriented framework for model management
  (e.g., \cite{muh93});
 ~{\it 5.} meta-modeling approaches (e.g., \cite{muh94});
 ~{\it 6.}  manipulation of composite models
  (e.g., \cite{ger90},\cite{lee93});
 ~{\it 7.} knowledge-based model management (e.g., \cite{lia88});
 and
 ~{\it 8.} applying machine learning to model management (e.g., \cite{sho88}).


 In the field of algorithms/algorithm systems
 the following interesting trends can be pointed out:

 {\it 1.} Automatic algorithm design
 as a combinatorial meta-problem  \cite{wallace02}.
 Here typical entries (for algorithm design) are:
 (i) choice of problem variables,
 (ii) choice of constraints, and
 (iii) choice of search method and constraint behavior.
 The approach realizes a joint design of problem
 formulation and algorithm.

 {\it 2.} Adaptive algorithms/software (e.g., \cite{mck04}).

 {\it 3.} Usage of generic library of problem solving methods/algorithms
  (e.g., \cite{kot12},\cite{raj06},\cite{rice76}).

 {\it 4.} Algorithm portfolio
 (e.g., \cite{gag10},\cite{gomes01},\cite{petrik06},\cite{zilber93}).

 {\it 5.} Reconfiguration or self-organizing algorithms
 (e.g., \cite{har01},\cite{zel04}).

 {\it 6.}  Cooperative/hybrid metaheuristics for combinatorial optimization
 (e.g., \cite{jour09},\cite{talbi09}).

 Some special research projects focus on the
 design of problem solving environments
 (e.g.,  \cite{dom02},\cite{fensel01},\cite{gal95}).
 For example,
 in \cite{raj06} a generic two-level library of problem solving
 is described which consists of two parts for the following:
 (a) basic task and subtasks problem solvers (task level),
 (b) problem solving technique (method level).
 Note special studies are needed for the problem formulation phase
 (e.g., \cite{eie03}).


 In the field of DSS,
 planning of decision making processes is a vital part
 of decision making engineering.
 In this case, design/building of  solving strategies
 is often under examination:
 ~(1) selection and integration of models from a model base
 (e.g.,
 \cite{dol84},\cite{gro91},\cite{kau94},\cite{liu88},\cite{vin91});
 ~(2) intelligent strategies for  decision making (e.g.,
 \cite{hanne01});
 ~(3)  MCDM techniques selection approaches
 (e.g., \cite{brau12},\cite{gui98},\cite{kornysh07},\cite{oz92});
  ~(4) expert-based hierarchical planning (e.g., \cite{vin91});
 ~(5) usage of decision making method families and their configuration
 (e.g., \cite{denkornysh09},\cite{denkornysh11},\cite{kornysh12});
 and
 ~(6) visual and interactive support for multicriteria
  decision making process.
  (e.g., \cite{trin05}).


 In this paper, an operational part of DSS for
 multicriteria ranking (sorting) and composition of corresponding
  composite solving strategies is examined.
 Usually, planning the decision support process
 is based on the selection and integration of models from a model base.
 These procedures use special model knowledge including descriptions
 of basic submodels and their connections, etc.
 (e.g., \cite{dol84}, \cite{kau94}, \cite{liu88}, \cite{vin91}).
 In general, four ``creative levels''can be considered:
     (i) a basic object,
     (ii) a selected object,
     (iii) a modified object, and
     (iv) a designed composite object (as composition of local elements/components).
 The levels are close to
  a ``scale of creativity''
 that has been  suggested by G.S. Altshuller
 (e.g., \cite{alt84})
 The levels may be applied to various parts of DSS
 (information, operations, user).
 Here, composition of solving strategies
  for multicriteria ranking (sorting)
 from some local techniques or procedures
   is studied.
 Our morphological composition approach is close to the strategy-generation table suggested
 by Howard \cite{how88}.
 Note here issues of creativity in decision making are not considered
 \cite{kee94}.

 The morphological approach for the design/planning of
 composite solving strategies is based on
 ``hierarchical morphological multicriteria design''
 (HMMD) which involves a series and/or parallel composition of
 the design alternatives (DAs) for data processing while
 taking into account interconnections (IC) among DAs
 (\cite{lev96a},\cite{lev98},\cite{lev06},\cite{lev12morph}).
 In the case of model management/engineering,
 DAs correspond to alternative models
 or algorithms/interactive procedures.
 HMMD is based on the following assumptions:
   (1) tree-like structure of a designed system;
   (2) system effectiveness (excellence) is represented as an
   aggregation of two parts:  effectiveness
   of subsystems (components),
   and effectiveness of compatibility among subsystems;
  (3) monotone criteria (Cr) for DAs;
   (4) ordinal coordinated priorities of DAs
  (\(\iota = \overline{1,l}\),  \(1\) corresponds to the best one);
  and
   (5) ordinal coordinated estimates of compatibility
   for each ``neighbor'' pair of DAs
  (\(w=\overline{1,\nu}\), \(\nu\) corresponds to the best one).
 In the case of interval multiset estimates
 (a generalization of ordinal estimates \cite{lev12a}),
 qualities of DAs and/or their compatibility
 are evaluated with usage of this kind of estimates
 (i.e., a special new type of poset-like scales).

  HMMD consists of the following main phases:


 {\bf Phase 1.} Top-Down design of tree-like system model
 including Cr and factors of compatibility.

 {\bf Phase 2.} Generation of DAs for leaf nodes of model.

 {\bf Phase 3.} Bottom-Up hierarchical selection and
 composition (iterative):
 (a) assessment of DAs on Cr,
 (b) assessment of IC,
 (c) computing the priorities of DAs, and
 (d) composing DAs for a higher hierarchical level.


 Note that similar approaches have been used in engineering design
 as combinatorial optimization models and morphological analysis
 (e.g., \cite{ayr69},\cite{ber93},\cite{gro90},\cite{jon92},\cite{zwi69}).

 Our study is based on the experience in the
 design and implementation
 of DSS COMBI for multicriteria ranking
  (joint project of the author and Andrew A. Michailov;
 1984...1991)
 (\cite{lev93},\cite{lev98},\cite{levmih88}).
 This DSS was based on a series-parallel transformation of
 preference relations
 (\cite{lev86},\cite{lev93},\cite{levmih88}).
 Operational environment of DSS COMBI includes a morphology of
 the composite solving strategy
 for forming/transforming the preference relations, linear orderings,
 and rankings.
 DSS COMBI corresponds to DSS generation 4.5
 of 5-level classification of Jelassi
 (\cite{jel86},\cite{sis99}).
 The following is described:
 (i) methodological issues in decision making technology for
 multicriteria ranking,
 (ii)  a scheme for designing a series solving strategy.

\section{METHODOLOGICAL ISSUES}

\subsection{Decision Making Framework, Solving Scheme, Problems}

 Generally,
 decision making process is based on the following basic parts
 (e.g., \cite{lev98},\cite{lev06},\cite{sim76}):
 (i) alternatives,
 (ii) criteria  and estimates of the alternatives upon the criteria,
 (iii) preferences over the alternatives,
 (iv) solving method(s),
 (v) decision(s), and
 (vi) specialists (decision maker, support experts).
 Fig. 1 depicts a basic framework of the decision making process.

\begin{center}
\begin{picture}(75,62)
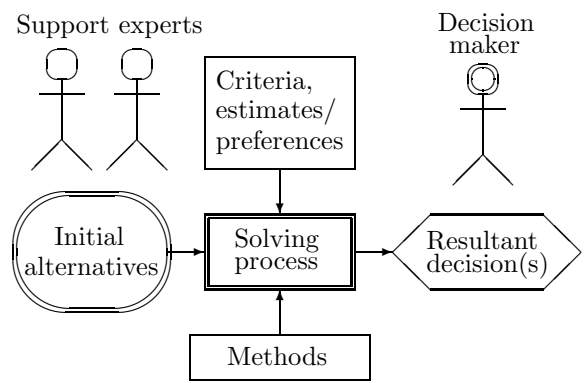
\put(03,00){\makebox(0,0)[bl] {Fig. 1. Basic framework of decision
making}}


\put(00,53.5){\makebox(0,0)[bl]{Support experts}}

\put(06,40){\line(0,1){8}}

\put(02,36){\line(1,1){4}} \put(10,36){\line(-1,1){4}}

\put(03,46){\line(1,0){6}}

\put(06,50){\oval(4,4)}

\put(16,40){\line(0,1){8}}

\put(12,36){\line(1,1){4}} \put(20,36){\line(-1,1){4}}

\put(13,46){\line(1,0){6}}

\put(16,50){\oval(4,4)}


\put(56,54.5){\makebox(0,0)[bl]{Decision}}
\put(57.5,51.5){\makebox(0,0)[bl]{maker}}

\put(62,38){\line(0,1){8}}

\put(58,34){\line(1,1){4}} \put(66,34){\line(-1,1){4}}

\put(59,44){\line(1,0){6}}

\put(62,48){\oval(4,4)} \put(62,48){\oval(3,3)}


\put(54.5,25.5){\makebox(0,0)[bl]{Resultant}}
\put(54.5,21.5){\makebox(0,0)[bl]{decision(s)}}

\put(50,25){\line(1,1){5}} \put(50,25){\line(1,-1){5}}

\put(75,25){\line(-1,1){5}} \put(75,25){\line(-1,-1){5}}

\put(55,30){\line(1,0){15}} \put(55,20){\line(1,0){15}}

\put(05,26){\makebox(0,0)[bl]{Initial}}
\put(01,22){\makebox(0,0)[bl]{alternatives}}

\put(10,25){\oval(21,16)} \put(10,25){\oval(20,15)}

\put(20,25){\vector(1,0){5}}

\put(26.5,46){\makebox(0,0)[bl]{Criteria,}}
\put(26.5,42){\makebox(0,0)[bl]{estimates/}}
\put(26.5,38){\makebox(0,0)[bl]{preferences}}

\put(25,36){\line(1,0){20}} \put(25,51){\line(1,0){20}}
\put(25,36){\line(0,1){15}} \put(45,36){\line(0,1){15}}
\put(35,36){\vector(0,-1){6}}

\put(29,25){\makebox(0,0)[bl]{Solving }}
\put(29,22){\makebox(0,0)[bl]{process}}
\put(25,20){\line(1,0){20}} \put(25,30){\line(1,0){20}}
\put(25,20){\line(0,1){10}} \put(45,20){\line(0,1){10}}
\put(25.5,20.5){\line(1,0){19}} \put(25.5,29.5){\line(1,0){19}}
\put(25.5,20.5){\line(0,1){9}} \put(44.5,20.5){\line(0,1){9}}
\put(45,25){\vector(1,0){5}}

\put(28,10){\makebox(0,0)[bl]{Methods}}
\put(23,08){\line(1,0){24}} \put(23,14){\line(1,0){24}}
\put(23,08){\line(0,1){6}} \put(47,08){\line(0,1){6}}
\put(35,14){\vector(0,1){6}}

\end{picture}
\end{center}

 H. Simon has suggested a rational decision making
 based on the choice problem  \cite{sim76}:
 (i) the identification and listing of  all the alternatives,
 (ii) determination of all the consequences resulting from each of the
 alternatives, and
 (iii) the comparison of the accuracy and efficiency of each of
 theses sets of consequences.
 A modified version of the approach is the following
 (\cite{lev86},\cite{lev06}):


 {\it Stage 1.} Analysis of the examined system/process,
 extraction of the problem.

 {\it Stage 2.} Problem structuring:
 (2.1.) generation of alternatives,
 (2.2.) generation of criteria and scale for each criterion.

 {\it Stage 3.} Obtaining the initial information
 (estimates of the alternatives, preferences over the
 alternatives).

 {\it Stage 4.} Solving process to obtain the decision(s).

 {\it Stage 5.} Analysis of the obtained decision(s).


 Fig. 2 depicts an extended decision making scheme including some
 feedback lines.

\begin{center}
\begin{picture}(82,143)
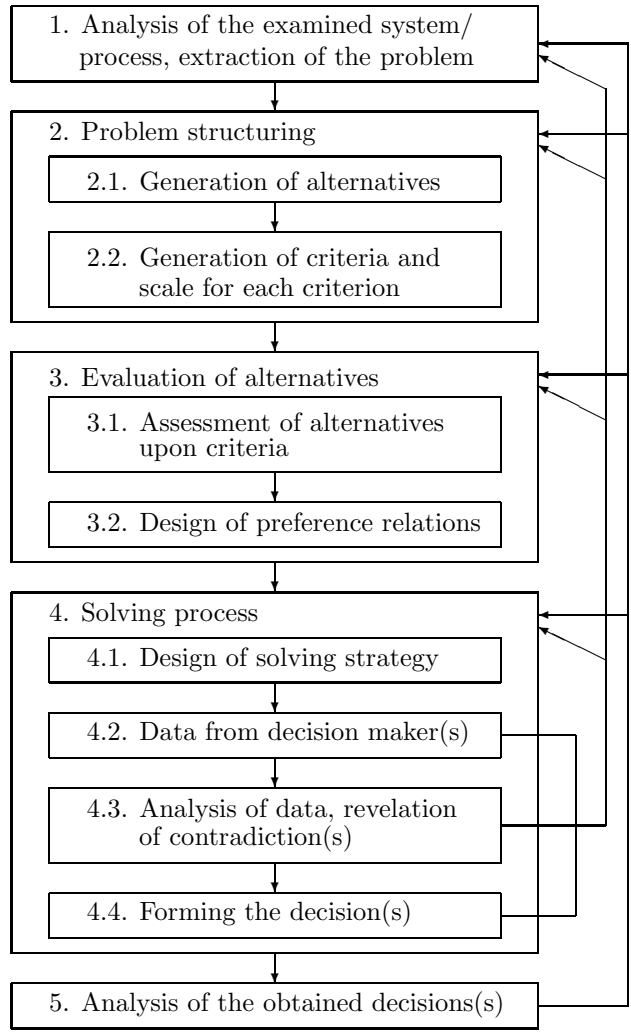
\put(07,00){\makebox(0,0)[bl] {Fig. 2. Extended decision making
scheme}}


\put(05,137.5){\makebox(0,0)[bl]{1. Analysis of the examined
 system/}}

\put(09,133.5){\makebox(0,0)[bl]{process, extraction of the
problem}}

\put(00,132){\line(1,0){70}} \put(00,142){\line(1,0){70}}
\put(00,132){\line(0,1){10}}  \put(70,132){\line(0,1){10}}


\put(05,123.5){\makebox(0,0)[bl]{2. Problem structuring}}

\put(00,100){\line(1,0){70}} \put(00,128){\line(1,0){70}}
\put(00,100){\line(0,1){28}}  \put(70,100){\line(0,1){28}}

\put(35,132){\vector(0,-1){4}}


\put(10,117.5){\makebox(0,0)[bl]{2.1. Generation of alternatives}}

\put(05,116){\line(1,0){60}} \put(05,122){\line(1,0){60}}
\put(05,116){\line(0,1){06}}  \put(65,116){\line(0,1){06}}


\put(10,107.5){\makebox(0,0)[bl]{2.2. Generation of criteria and}}

\put(17,103.5){\makebox(0,0)[bl]{scale for each criterion}}

\put(05,102){\line(1,0){60}} \put(05,112){\line(1,0){60}}
\put(05,102){\line(0,1){10}}  \put(65,102){\line(0,1){10}}

\put(35,116){\vector(0,-1){4}}


\put(05,91.5){\makebox(0,0)[bl]{3. Evaluation of alternatives}}

\put(00,68){\line(1,0){70}} \put(00,96){\line(1,0){70}}
\put(00,68){\line(0,1){28}}  \put(70,68){\line(0,1){28}}

\put(35,100){\vector(0,-1){4}}


\put(10,85.5){\makebox(0,0)[bl]{3.1. Assessment of alternatives}}

\put(17,81.5){\makebox(0,0)[bl]{upon criteria}}

\put(05,80){\line(1,0){60}} \put(05,90){\line(1,0){60}}
\put(05,80){\line(0,1){10}}  \put(65,80){\line(0,1){10}}


\put(10,71.5){\makebox(0,0)[bl]{3.2. Design of preference
relations}}

\put(05,70){\line(1,0){60}} \put(05,76){\line(1,0){60}}
\put(05,70){\line(0,1){06}}  \put(65,70){\line(0,1){06}}

\put(35,80){\vector(0,-1){4}}


\put(05,59.5){\makebox(0,0)[bl]{4. Solving process}}

\put(00,16){\line(1,0){70}} \put(00,64){\line(1,0){70}}
\put(00,16){\line(0,1){48}}  \put(70,16){\line(0,1){48}}

\put(35,68){\vector(0,-1){4}}


\put(10,53.5){\makebox(0,0)[bl]{4.1. Design of solving strategy}}

\put(05,52){\line(1,0){60}} \put(05,58){\line(1,0){60}}
\put(05,52){\line(0,1){06}}  \put(65,52){\line(0,1){06}}



\put(10,43.5){\makebox(0,0)[bl]{4.2. Data from decision maker(s)}}

\put(05,42){\line(1,0){60}} \put(05,48){\line(1,0){60}}
\put(05,42){\line(0,1){06}}  \put(65,42){\line(0,1){06}}

\put(35,52){\vector(0,-1){4}}


\put(10,33.5){\makebox(0,0)[bl]{4.3. Analysis of data,
revelation}}

\put(17,29.5){\makebox(0,0)[bl]{of contradiction(s)}}

\put(05,28){\line(1,0){60}} \put(05,38){\line(1,0){60}}
\put(05,28){\line(0,1){10}}  \put(65,28){\line(0,1){10}}

\put(35,42){\vector(0,-1){4}}


\put(10,19.5){\makebox(0,0)[bl]{4.4. Forming the decision(s)}}

\put(05,18){\line(1,0){60}} \put(05,24){\line(1,0){60}}
\put(05,18){\line(0,1){6}}  \put(65,18){\line(0,1){6}}

\put(35,28){\vector(0,-1){4}}


\put(05,07.5){\makebox(0,0)[bl]{5. Analysis of the obtained
decisions(s)}}

\put(00,06){\line(1,0){70}} \put(00,12){\line(1,0){70}}
\put(00,06){\line(0,1){06}}  \put(70,06){\line(0,1){06}}


\put(35,16){\vector(0,-1){4}}


\put(70,09){\line(1,0){12}} \put(82,09){\line(0,1){128}}

\put(82,137){\vector(-1,0){12}}

\put(82,125){\vector(-1,0){12}}

\put(82,93){\vector(-1,0){12}}

\put(82,61){\vector(-1,0){12}}


\put(65,45){\line(1,0){10}} \put(65,33){\line(1,0){14}}
\put(65,21){\line(1,0){10}} \put(75,21){\line(0,1){24}}

\put(79,33){\line(0,1){98}}

\put(79,131){\vector(-2,1){9}}

\put(79,119){\vector(-2,1){9}}

\put(79,87){\vector(-2,1){9}}

\put(79,55){\vector(-2,1){9}}

\end{picture}
\end{center}

 Fig. 3 depicts the basic decision making problems
 (\cite{fig05},\cite{kee76},\cite{lev86},\cite{lev98},\cite{mirkin79},\cite{roy96}):

 (a) choice/selection
 (e.g.,
 \cite{aiz95},\cite{fig05},\cite{kee76},\cite{mirkin79},\cite{sim76},\cite{tver72}),

 (b) linear ordering
 (e.g.,
 \cite{fig05},\cite{fis70},\cite{kee76},\cite{mirkin79}),

 (c) sorting/ranking
 (e.g., \cite{bregar08},\cite{chen07},\cite{fig05},\cite{lev98},\cite{levmih88},\cite{mirkin79},\cite{roy96},\cite{zop00},\cite{zop02}),
and

 (d) clustering/classification
 (e.g.,
 \cite{chen06},\cite{jain99},\cite{mirkin79},\cite{mirkin05},\cite{zop02}).

 In general, it is reasonable to extend the decision making
 process by additional management/monitoring
 and support analysis/learning of user(s)
 (Fig. 4).
 This approach  was implemented in DSS COMBI
  (\cite{lev86},\cite{lev93},\cite{lev98},\cite{levmih87a},\cite{levmih87b},\cite{levmih87c},\cite{levmih88},\cite{levmih90}),
  for example:
 (i) special hypertext system
 for learning and support,
 (ii) analysis and diagnosis of user(s),
 (iii) library of typical DM problems for various domains,
 (iv) basic typical solving composite strategies, and
 (v) possibility for retrieval and analysis of intermediate
 information/solutions.

 \begin{center}
\begin{picture}(110,74)
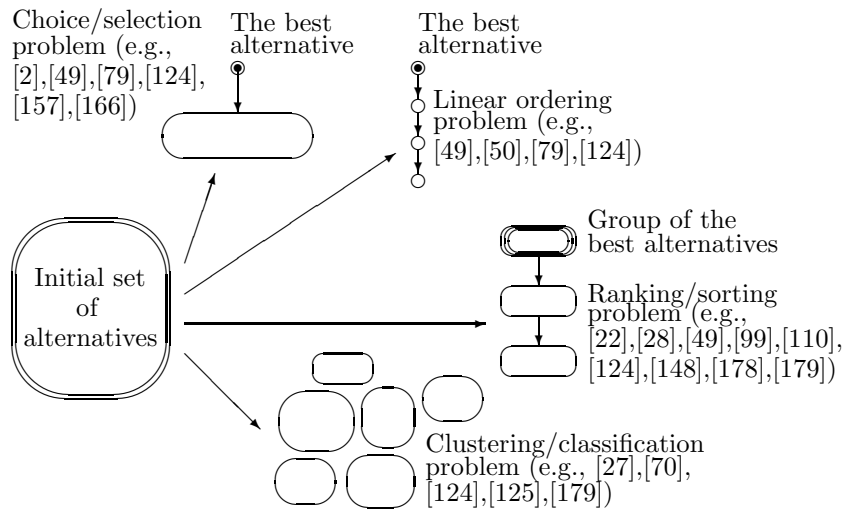
\put(18,0){\makebox(0,0)[bl]{Fig. 3. Basic decision making
problems}}

\put(10.5,32){\oval(21,24)} \put(10.5,32){\oval(20,23)}

\put(3,35){\makebox(0,0)[bl]{Initial set}}
\put(1.5,31){\makebox(0,0)[bl]{~~~~~ of }}
\put(1.5,27){\makebox(0,0)[bl]{alternatives}}

\put(23,38){\vector(1,3){04}}

\put(23,34){\vector(3,2){28}}

\put(23,30){\vector(1,0){40}}

\put(23,26){\vector(1,-1){10}}


\put(00,68.5){\makebox(0,0)[bl]{Choice/selection}}
\put(00,65){\makebox(0,0)[bl]{problem (e.g., }}
\put(00,61){\makebox(0,0)[bl]{\cite{aiz95},\cite{fig05},\cite{kee76},\cite{mirkin79},}}
\put(00,57){\makebox(0,0)[bl]{\cite{sim76},\cite{tver72})}}

\put(30,64){\circle{1.8}} \put(30,64){\circle*{1.0}}

\put(29,69){\makebox(0,0)[bl]{The best }}
\put(29,66){\makebox(0,0)[bl]{alternative}}

\put(30,64){\vector(0,-1){6}} \put(30,55){\oval(20,6)}


\put(55,12){\makebox(0,0)[bl]{Clustering/classification}}
\put(55,09){\makebox(0,0)[bl]{problem (e.g.,
 \cite{chen06},\cite{jain99},}}

\put(55,05.5){\makebox(0,0)[bl]{\cite{mirkin79},\cite{mirkin05},\cite{zop02})}}

\put(39,09){\oval(08,6)} \put(49,09){\oval(09,7)}

\put(40.5,17){\oval(10,8)} \put(50,17.5){\oval(07,8)}
\put(58.6,20){\oval(08,6)}

\put(44,24){\oval(08,4)}


\put(76.5,42){\makebox(0,0)[bl]{Group of the}}
\put(76.5,39.5){\makebox(0,0)[bl]{best alternatives}}

\put(70,41){\oval(09,3.5)}  \put(70,41){\oval(08,3)}

\put(70,41){\oval(10,4)} \put(70,39){\vector(0,-1){4}}
\put(70,33){\oval(10,4)} \put(70,31){\vector(0,-1){4}}
\put(70,25){\oval(10,4)}

\put(76.5,32){\makebox(0,0)[bl]{Ranking/sorting}}
\put(76.5,29.5){\makebox(0,0)[bl]{problem (e.g., }}
\put(76.5,26){\makebox(0,0)[bl]{\cite{bregar08},\cite{chen07},\cite{fig05},\cite{lev98},\cite{levmih88},}}
\put(76.5,22){\makebox(0,0)[bl]{\cite{mirkin79},\cite{roy96},\cite{zop00},\cite{zop02})}}


\put(56,58){\makebox(0,0)[bl]{Linear ordering }}
\put(56,55){\makebox(0,0)[bl]{problem (e.g., }}
\put(56,51){\makebox(0,0)[bl]{\cite{fig05},\cite{fis70},\cite{kee76},\cite{mirkin79})}}

\put(54,69){\makebox(0,0)[bl]{The best }}
\put(54,66){\makebox(0,0)[bl]{alternative}}

\put(54,64){\circle{1.8}} \put(54,64){\circle*{1.0}}
\put(54,63){\vector(0,-1){3}}

\put(54,59){\circle{1.8}} \put(54,58){\vector(0,-1){3}}

\put(54,54){\circle{1.8}} \put(54,53){\vector(0,-1){3}}

\put(54,49){\circle{1.8}}

\end{picture}
\end{center}

\begin{center}
\begin{picture}(124,90)
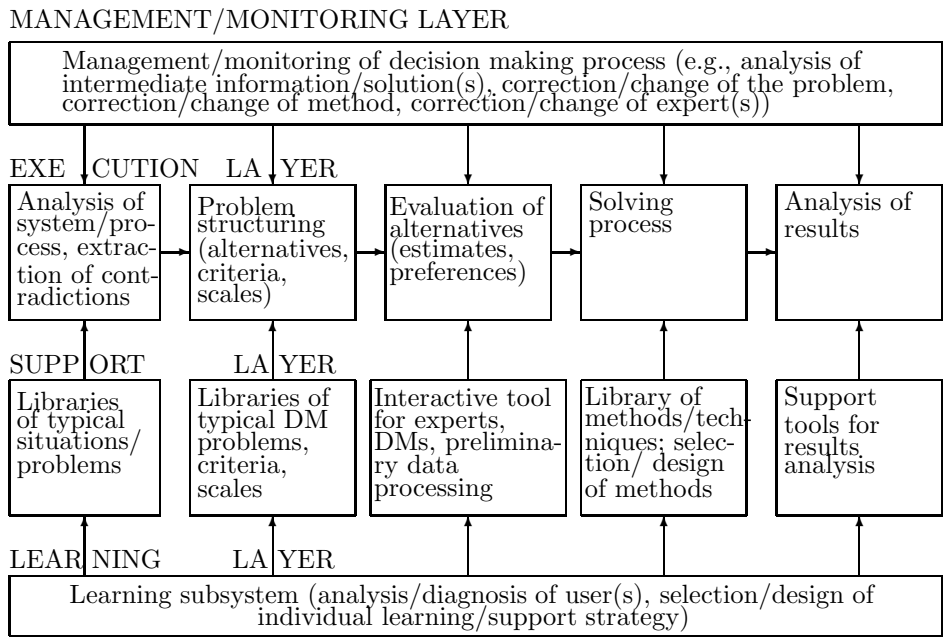
\put(09.5,00){\makebox(0,0)[bl] {Fig. 4. Extended four-layer
 architecture of decision making process}}


\put(00,66){\makebox(0,0)[bl]{EXE~~~CUTION ~~LA~~YER}}

\put(01,61){\makebox(0,0)[bl]{Analysis of}}
\put(01,58){\makebox(0,0)[bl]{system/pro-}}
\put(01,55){\makebox(0,0)[bl]{cess, extrac-}}
\put(01,52){\makebox(0,0)[bl]{tion of cont-}}
\put(01,49){\makebox(0,0)[bl]{radictions}}

\put(00,47){\line(1,0){20}} \put(00,65){\line(1,0){20}}
\put(00,47){\line(0,1){18}} \put(20,47){\line(0,1){18}}

\put(20,56){\vector(1,0){4}}


\put(25,61){\makebox(0,0)[bl]{Problem }}
\put(25,58){\makebox(0,0)[bl]{structuring}}
\put(25,54.5){\makebox(0,0)[bl]{(alternatives,}}
\put(25,51.5){\makebox(0,0)[bl]{criteria,}}
\put(25,48.5){\makebox(0,0)[bl]{scales)}}

\put(24,47){\line(1,0){22}} \put(24,65){\line(1,0){22}}
\put(24,47){\line(0,1){18}} \put(46,47){\line(0,1){18}}

\put(46,56){\vector(1,0){4}}


\put(50.5,61){\makebox(0,0)[bl]{Evaluation of}}
\put(50.5,58){\makebox(0,0)[bl]{alternatives}}
\put(50.5,54.5){\makebox(0,0)[bl]{(estimates,}}
\put(50.5,51.5){\makebox(0,0)[bl]{preferences)}}

\put(50,47){\line(1,0){22}} \put(50,65){\line(1,0){22}}
\put(50,47){\line(0,1){18}} \put(72,47){\line(0,1){18}}

\put(72,56){\vector(1,0){4}}


\put(77,61){\makebox(0,0)[bl]{Solving}}
\put(77,58){\makebox(0,0)[bl]{process}}
\put(77,55){\makebox(0,0)[bl]{}} \put(77,52){\makebox(0,0)[bl]{}}

\put(76,47){\line(1,0){22}} \put(76,65){\line(1,0){22}}
\put(76,47){\line(0,1){18}} \put(98,47){\line(0,1){18}}

\put(98,56){\vector(1,0){4}}


\put(103,61){\makebox(0,0)[bl]{Analysis of }}
\put(103,58){\makebox(0,0)[bl]{results}}
\put(103,55){\makebox(0,0)[bl]{}} \put(93,52){\makebox(0,0)[bl]{}}

\put(102,47){\line(1,0){22}} \put(102,65){\line(1,0){22}}
\put(102,47){\line(0,1){18}} \put(124,47){\line(0,1){18}}


\put(00,40){\makebox(0,0)[bl]{SUPP~ORT ~~~~~~~~ LA~YER}}

\put(01,35){\makebox(0,0)[bl]{Libraries }}
\put(01,32){\makebox(0,0)[bl]{of typical }}
\put(01,29){\makebox(0,0)[bl]{situations/}}
\put(01,26){\makebox(0,0)[bl]{problems}}
\put(01,23){\makebox(0,0)[bl]{}}

\put(00,21){\line(1,0){20}} \put(00,39){\line(1,0){20}}
\put(00,21){\line(0,1){18}} \put(20,21){\line(0,1){18}}

\put(10,39){\vector(0,1){8}}


\put(25,35.5){\makebox(0,0)[bl]{Libraries of}}
\put(25,32){\makebox(0,0)[bl]{typical DM}}
\put(25,29){\makebox(0,0)[bl]{problems, }}
\put(25,26){\makebox(0,0)[bl]{criteria,}}
\put(25,23.5){\makebox(0,0)[bl]{scales }}

\put(24,21){\line(1,0){22}} \put(24,39){\line(1,0){22}}
\put(24,21){\line(0,1){18}} \put(46,21){\line(0,1){18}}

\put(35,39){\vector(0,1){8}}


\put(48.5,35.5){\makebox(0,0)[bl]{Interactive tool}}
\put(48.5,32){\makebox(0,0)[bl]{for experts, }}
\put(48.5,29){\makebox(0,0)[bl]{DMs, prelimina-}}
\put(48.5,26){\makebox(0,0)[bl]{ry data}}
\put(48.5,23){\makebox(0,0)[bl]{processing}}

\put(48,21){\line(1,0){26}} \put(48,39){\line(1,0){26}}
\put(48,21){\line(0,1){18}} \put(74,21){\line(0,1){18}}

\put(61,39){\vector(0,1){8}}


\put(76.5,35){\makebox(0,0)[bl]{Library of}}
\put(76.5,32){\makebox(0,0)[bl]{methods/tech- }}
\put(76.5,29){\makebox(0,0)[bl]{niques; selec-}}
\put(76.5,26){\makebox(0,0)[bl]{tion/ design}}
\put(76.5,23.5){\makebox(0,0)[bl]{of methods}}

\put(76,21){\line(1,0){22}} \put(76,39){\line(1,0){22}}
\put(76,21){\line(0,1){18}} \put(98,21){\line(0,1){18}}

\put(87,39){\vector(0,1){8}}


\put(103,35){\makebox(0,0)[bl]{Support}}
\put(103,32){\makebox(0,0)[bl]{tools for}}
\put(103,29){\makebox(0,0)[bl]{results}}
\put(103,26){\makebox(0,0)[bl]{analysis}}

\put(102,21){\line(1,0){22}} \put(102,39){\line(1,0){22}}
\put(102,21){\line(0,1){18}} \put(124,21){\line(0,1){18}}

\put(113,39){\vector(0,1){8}}


\put(00,14){\makebox(0,0)[bl]{LEAR~NING ~~~~~~ LA~YER}}

\put(08,8.5){\makebox(0,0)[bl]{Learning subsystem
 (analysis/diagnosis of user(s), selection/design  of }}

\put(33,05.5){\makebox(0,0)[bl]{individual learning/support
 strategy) }}

\put(00,05){\line(1,0){124}} \put(00,13){\line(1,0){124}}
\put(00,05){\line(0,1){8}} \put(124,05){\line(0,1){8}}

\put(10,13){\vector(0,1){8}} \put(35,13){\vector(0,1){8}}
\put(61,13){\vector(0,1){8}} \put(87,13){\vector(0,1){8}}
\put(113,13){\vector(0,1){8}}


\put(00,85){\makebox(0,0)[bl]{MANAGEMENT/MONITORING LAYER}}

\put(07,79.5){\makebox(0,0)[bl]{Management/monitoring of decision
 making process (e.g., analysis of  }}

\put(07,76.5){\makebox(0,0)[bl]{intermediate
 information/solution(s),
 correction/change of the problem, }}

\put(07,74){\makebox(0,0)[bl]{correction/change of method,
correction/change of expert(s)) }}

\put(00,073){\line(1,0){124}} \put(00,84){\line(1,0){124}}
\put(00,073){\line(0,1){11}} \put(124,073){\line(0,1){11}}

\put(10,73){\vector(0,-1){8}} \put(35,73){\vector(0,-1){8}}
\put(61,73){\vector(0,-1){8}} \put(87,73){\vector(0,-1){8}}
\put(113,73){\vector(0,-1){8}}

\end{picture}
\end{center}

\subsection{Multicriteria Ranking (Sorting)}

 Let \(A=\{1,...,i,...,n\}\) be a set of alternatives
 (items)
 which are evaluated on
 criteria
 \( K=  \{ 1,...,j,...,d \}\), and
 \(z(i,j)\) is an estimate (quantitative, ordinal) of alternative \(i\)
 on criterion \(j\).
 The matrix \(Z = \{z(i,j)\}\) may be
 mapped into a poset on \(A\).
 We search for the following resultant kind of the poset as a partition
 with ordered subsets (a layered structure):
   \(B=\{ B_{1},...,B_{k},...,B_{m}\} \),
   \( B_{k_{1}}\& B_{k_{2}}=\O \) if \(k_{1}\neq k_{2}\),
   and each alternative from \(B_{k_{1}}\) (layer \(k_{1}\))
 dominates each alternative from \(B_{k_{2}}\) (layer \(k_{2}\)),
 if \(k_{1}\leq k_{2}\).
   Thus, each alternative has a priority which equals the number of the
 corresponding layer.
 This problem is illustrated in Fig. 5.
 This problem belongs to a class of ill-structured problems
 by classification of H. Simon \cite{simon58}.
 In general, the resultant ordered subsets can have intersections
 (i.e., the problem can be targeted to obtain interval priorities
 for the alternatives) (Fig. 6).
 In this case, the ``fuzzy'' decision will be denoted by
 \(\widetilde{B}\)
 (a layered structure with intersection of the layers).
 Let \(\overline{B}\) be a linear ordering of alternatives.

 The basic techniques for the multicriteria ranking
 (sorting) problems are
 the following (e.g., \cite{bue92},\cite{dye92},\cite{zav03}):
  (1) statistical approaches \cite{boy74};
  (2) multi-attribute utility analysis (e.g., \cite{fis70},\cite{kee76});
  (3) multi-criterion decision making (e.g., \cite{kor84});
  (4) analytic hierarchy process (e.g., \cite{saa80});
  (5) outranking techniques (e.g.,  \cite{beh10},\cite{roy96});
  (6) knowledge bases (e.g., \cite{li87});
  (7) neural network (e.g., \cite{wan92});
 (8) logical procedures (e.g., \cite{levmih87c});
  (9) expert judgment (e.g., \cite{larmf86}); and
  (10) hybrid techniques (e.g., \cite{fed93}).
 In the main, the above-mentioned techniques correspond
 to one-phase problem solving framework.

\begin{center}
\begin{picture}(72,45)

\put(00,00){\makebox(0,0)[bl]{Fig. 5. Multicriteria ranking
(sorting)}}

\put(14.5,25){\oval(29,30)}

\put(5,28){\makebox(0,0)[bl]{Initial set of}}
\put(5,24){\makebox(0,0)[bl]{alternatives}}
\put(0.5,18.5){\makebox(0,0)[bl]{\(A=\{1,...,i,...,n\}\)}}

\put(31,32){\vector(1,1){6}}

\put(31,25){\vector(1,0){6}}

\put(31,18){\vector(1,-1){6}}


\put(51,42){\oval(20,6)}
\put(49,40.5){\makebox(0,0)[bl]{\(B_{1}\)}}

\put(48,39){\vector(0,-1){4}} \put(51,39){\vector(0,-1){4}}
\put(54,39){\vector(0,-1){4}}


\put(47.5,33){\makebox(0,0)[bl]{{\bf . . .}}}

\put(48,32){\vector(0,-1){4}} \put(51,32){\vector(0,-1){4}}
\put(54,32){\vector(0,-1){4}}

\put(51,25){\oval(20,6)}

\put(49,23.5){\makebox(0,0)[bl]{\(B_{k}\)}}

\put(48,22){\vector(0,-1){4}} \put(51,22){\vector(0,-1){4}}
\put(54,22){\vector(0,-1){4}}


\put(47.5,16){\makebox(0,0)[bl]{{\bf . . .}}}

\put(48,15){\vector(0,-1){4}} \put(51,15){\vector(0,-1){4}}
\put(54,15){\vector(0,-1){4}}

\put(51,8){\oval(20,6)}

\put(49,6.5){\makebox(0,0)[bl]{\(B_{m}\)}}

\end{picture}
%
\begin{picture}(62,45)

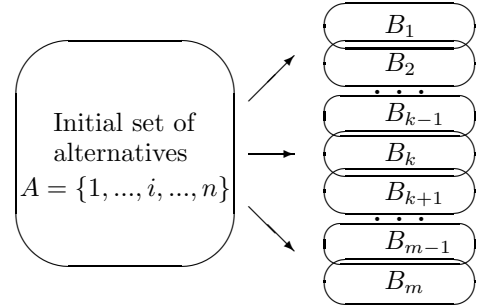
\put(00,00){\makebox(0,0)[bl]{Fig. 6. Ranking with interval
priorities}}

\put(14.5,25){\oval(29,30)}

\put(5,28){\makebox(0,0)[bl]{Initial set of}}
\put(5,24){\makebox(0,0)[bl]{alternatives}}
\put(0.5,18.5){\makebox(0,0)[bl]{\(A=\{1,...,i,...,n\}\)}}


\put(31,32){\vector(1,1){6}}

\put(31,25){\vector(1,0){6}}

\put(31,18){\vector(1,-1){6}}


\put(51,42){\oval(20,6)}
\put(49,40.5){\makebox(0,0)[bl]{\(B_{1}\)}}

\put(51,37){\oval(20,6)}
\put(49,35.5){\makebox(0,0)[bl]{\(B_{2}\)}}


\put(47.5,33){\makebox(0,0)[bl]{{\bf . . .}}}


\put(51,30){\oval(20,5.4)}

\put(49,28.5){\makebox(0,0)[bl]{\(B_{k-1}\)}}


\put(51,25){\oval(20,6)}

\put(49,23.5){\makebox(0,0)[bl]{\(B_{k}\)}}


\put(51,20){\oval(20,6)}

\put(49,18){\makebox(0,0)[bl]{\(B_{k+1}\)}}


\put(47.5,16){\makebox(0,0)[bl]{{\bf . . .}}}

\put(51,13){\oval(20,5.4)}

\put(49,11.5){\makebox(0,0)[bl]{\(B_{m-1}\)}}

\put(51,8){\oval(20,6)}

\put(49,6.5){\makebox(0,0)[bl]{\(B_{m}\)}}

\end{picture}
\end{center}

 The following numerical example illustrates the examined
 multicriteria problem.

 Let
 \(A = \{ A_{1},A_{2},A_{3},A_{4},A_{5},A_{6},A_{7},A_{8},A_{9}
 \}\) be a set of initial alternatives.
 The alternatives are evaluated upon 2 criteria
  \(K = \{ K_{1},K_{2} \}\),
 an ordinal scale is used for each criterion
 (\([0,1,2,3,4]\) and \(4\) corresponds to the best level).
 The ordinal estimates of alternatives are presented in Table 1,
 a space of estimates is depicted in Fig. 7.
 Three types of solutions are presented:

 (a) linear ordering  \(\overline{B}\)
  (Fig. 8;
 e.g., additive utility function is used);

 (b) ranking (sorting, four linear ordered subsets)
 \(B\)
 (Fig. 9; e.g., expert judgement is used);
 and

 (c) ``fuzzy'' ranking
 \(\widetilde{B}\)
 (Fig. 10; e.g., expert judgement is used).

 Table 2 integrates the obtained priorities of alternatives.

\begin{center}
\begin{picture}(54,57)

\put(05.5,52){\makebox(0,0)[bl]{Table 1. Estimates}}

\put(00,00){\line(1,0){42}} \put(00,38){\line(1,0){42}}
\put(10,44){\line(1,0){32}} \put(00,50){\line(1,0){42}}

\put(00,0){\line(0,1){50}} \put(10,0){\line(0,1){50}}
\put(42,0){\line(0,1){50}}

\put(01,34){\makebox(0,0)[bl]{\(A_{1}\)}}
\put(01,30){\makebox(0,0)[bl]{\(A_{2}\)}}
\put(01,26){\makebox(0,0)[bl]{\(A_{3}\)}}
\put(01,22){\makebox(0,0)[bl]{\(A_{4}\)}}
\put(01,18){\makebox(0,0)[bl]{\(A_{5}\)}}

\put(01,14){\makebox(0,0)[bl]{\(A_{6}\)}}
\put(01,10){\makebox(0,0)[bl]{\(A_{7}\)}}
\put(01,06){\makebox(0,0)[bl]{\(A_{8}\)}}
\put(01,02){\makebox(0,0)[bl]{\(A_{9}\)}}

\put(26,38){\line(0,1){6}}

\put(0.5,46){\makebox(0,0)[bl]{Alter-}}
\put(0.5,42){\makebox(0,0)[bl]{native}}

\put(20,46){\makebox(0,0)[bl]{Criteria}}

\put(16,40){\makebox(0,0)[bl]{\(K_{1}\)}}

\put(32,40){\makebox(0,0)[bl]{\(K_{2}\)}}

\put(17,34){\makebox(0,0)[bl]{\(2\)}}
\put(33,34){\makebox(0,0)[bl]{\(3\)}}

\put(17,30){\makebox(0,0)[bl]{\(2\)}}
\put(33,30){\makebox(0,0)[bl]{\(4\)}}

\put(17,26){\makebox(0,0)[bl]{\(1\)}}
\put(33,26){\makebox(0,0)[bl]{\(3\)}}

\put(17,22){\makebox(0,0)[bl]{\(4\)}}
\put(33,22){\makebox(0,0)[bl]{\(4\)}}

\put(17,18){\makebox(0,0)[bl]{\(1\)}}
\put(33,18){\makebox(0,0)[bl]{\(1\)}}

\put(17,14){\makebox(0,0)[bl]{\(4\)}}
\put(33,14){\makebox(0,0)[bl]{\(3\)}}

\put(17,10){\makebox(0,0)[bl]{\(2\)}}
\put(33,10){\makebox(0,0)[bl]{\(2\)}}

\put(17,06){\makebox(0,0)[bl]{\(0\)}}
\put(33,06){\makebox(0,0)[bl]{\(2\)}}

\put(17,02){\makebox(0,0)[bl]{\(2\)}}
\put(33,02){\makebox(0,0)[bl]{\(1\)}}

\end{picture}
%
\begin{picture}(57,57)

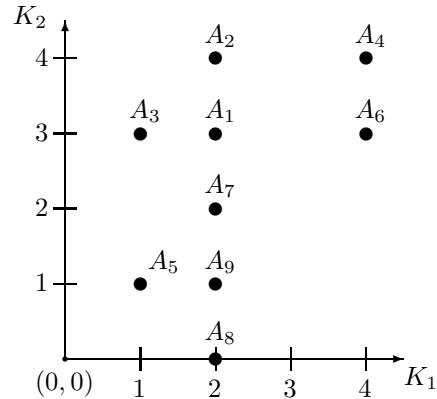
\put(010,00){\makebox(0,0)[bl]{Fig 7. Space of estimates}}

\put(003,54){\makebox(0,0)[bl]{\(K_{2}\)}}

\put(55,06){\makebox(0,0)[bl]{\(K_{1}\)}}

\put(06,05){\makebox(0,0)[bl]{\((0,0)\)}}

\put(10,10){\circle*{0.8}}

\put(10,10){\vector(0,1){45}} \put(10,10){\vector(1,0){45}}

\put(06,19){\makebox(0,0)[bl]{\(1\)}} \put(8.5,20){\line(1,0){3}}
\put(06,29){\makebox(0,0)[bl]{\(2\)}} \put(8.5,30){\line(1,0){3}}
\put(06,39){\makebox(0,0)[bl]{\(3\)}} \put(8.5,40){\line(1,0){3}}
\put(06,49){\makebox(0,0)[bl]{\(4\)}} \put(8.5,50){\line(1,0){3}}

\put(19,05){\makebox(0,0)[bl]{\(1\)}} \put(20,8.5){\line(0,1){3}}
\put(29,05){\makebox(0,0)[bl]{\(2\)}} \put(30,8.5){\line(0,1){3}}
\put(39,05){\makebox(0,0)[bl]{\(3\)}} \put(40,8.5){\line(0,1){3}}
\put(49,05){\makebox(0,0)[bl]{\(4\)}} \put(50,8.5){\line(0,1){3}}


\put(28.5,51.5){\makebox(0,0)[bl]{\(A_{2}\)}}
\put(30,50){\circle*{1.8}}


\put(48.5,51.5){\makebox(0,0)[bl]{\(A_{4}\)}}
\put(50,50){\circle*{1.8}}


\put(48.5,41.5){\makebox(0,0)[bl]{\(A_{6}\)}}
\put(50,40){\circle*{1.8}}



\put(18.5,41.5){\makebox(0,0)[bl]{\(A_{3}\)}}
\put(20,40){\circle*{1.8}}


\put(28.5,41.5){\makebox(0,0)[bl]{\(A_{1}\)}}
\put(30,40){\circle*{1.8}}


\put(28.5,31.5){\makebox(0,0)[bl]{\(A_{7}\)}}
\put(30,30){\circle*{1.8}}


\put(21,21.5){\makebox(0,0)[bl]{\(A_{5}\)}}
\put(20,20){\circle*{1.8}}


\put(28.5,12){\makebox(0,0)[bl]{\(A_{8}\)}}
\put(30,10){\circle*{1.8}}


\put(28.5,21.5){\makebox(0,0)[bl]{\(A_{9}\)}}
\put(30,20){\circle*{1.8}}

\end{picture}
\end{center}

\begin{center}
\begin{picture}(50,62)

\put(00,00){\makebox(0,0)[bl]{Fig 8. Linear ordering
 \(\overline{B}\)}}


\put(12,56.4){\makebox(0,0)[bl]{\(A_{4}\)}}

\put(23,58){\circle*{1.8}}


\put(23,58){\vector(0,-1){5}}

\put(12,50.4){\makebox(0,0)[bl]{\(A_{6}\)}}

\put(23,52){\circle*{1.8}}


\put(23,52){\vector(0,-1){5}}

\put(12,44.4){\makebox(0,0)[bl]{\(A_{2}\)}}

\put(23,46){\circle*{1.8}}


\put(23,46){\vector(0,-1){5}}

\put(12,38.4){\makebox(0,0)[bl]{\(A_{1}\)}}

\put(23,40){\circle*{1.8}}


\put(23,40){\vector(0,-1){5}}

\put(12,32.4){\makebox(0,0)[bl]{\(A_{3}\)}}

\put(23,34){\circle*{1.8}}


\put(23,34){\vector(0,-1){5}}

\put(12,26.4){\makebox(0,0)[bl]{\(A_{7}\)}}

\put(23,28){\circle*{1.8}}


\put(23,28){\vector(0,-1){5}}

\put(12,20.4){\makebox(0,0)[bl]{\(A_{9}\)}}

\put(23,22){\circle*{1.8}}


\put(23,22){\vector(0,-1){5}}

\put(12,14.4){\makebox(0,0)[bl]{\(A_{5}\)}}

\put(23,16){\circle*{1.8}}


\put(23,16){\vector(0,-1){5}}

\put(12,08.4){\makebox(0,0)[bl]{\(A_{8}\)}}

\put(23,10){\circle*{1.8}}


\end{picture}
%
\begin{picture}(50,50)

\put(12,00){\makebox(0,0)[bl]{Fig 9. Ranking \(B\)}}


\put(00,38.5){\makebox(0,0)[bl]{\(B_{1}:\)}}

\put(29,38.5){\makebox(0,0)[bl]{\(A_{4}\)}}

\put(30,40){\oval(40,6)}


\put(00,28.5){\makebox(0,0)[bl]{\(B_{2}:\)}}

\put(30,37){\vector(0,-1){4}}

\put(25,28.5){\makebox(0,0)[bl]{\(A_{2},A_{6}\)}}

\put(30,30){\oval(40,6)}


\put(00,18.5){\makebox(0,0)[bl]{\(B_{3}:\)}}

\put(30,27){\vector(0,-1){4}}

\put(23,18.5){\makebox(0,0)[bl]{\(A_{1},A_{3},A_{7}\)}}

\put(30,20){\oval(40,6)}


\put(00,08.5){\makebox(0,0)[bl]{\(B_{4}:\)}}

\put(30,17){\vector(0,-1){4}}

\put(23,08.5){\makebox(0,0)[bl]{\(A_{5},A_{8},A_{9}\)}}

\put(30,10){\oval(40,6)}

\end{picture}
\end{center}

 An approach for decision making
 as series or series-parallel processing (transformation, articulation)
 of preferences
 (combinatorial models for decision making)
 for suggested  in \cite{lev84}.
 In DSS COMBI,
 a functional graph has been suggested
 as impementation of the approach
 for multicriteria ranking
 (\cite{lev93},\cite{lev98},\cite{levmih88}).
 The graph was realized as a graphical menu.
 In this case, the solving strategy is
  combined  from a set of basic operations
  (local techniques or procedures)
 (e.g., forming preference relations over the alternatives,
 forming the intermediate linear ordering of the alternatives,
 forming the resultant decision structure over the alternatives)
 (Fig. 11).

\begin{center}
\begin{picture}(60,50)

\put(05,00){\makebox(0,0)[bl]{Fig 10. ``Fuzzy'' ranking
 \(\widetilde{B}\)}}


\put(00,40.5){\makebox(0,0)[bl]{\(\widetilde{B}_{1}:\)}}

\put(19,40.5){\makebox(0,0)[bl]{\(A_{4}\)}}

\put(30,42){\oval(40,6)}


\put(39.5,36){\oval(21,09)}

\put(30.5,34.5){\makebox(0,0)[bl]{\(\widetilde{B}_{1}\&\widetilde{B}_{2}:
A_{6}\)}}


\put(00,29.5){\makebox(0,0)[bl]{\(\widetilde{B}_{2}:\)}}

\put(20,39){\vector(0,-1){5}}

\put(19,29.5){\makebox(0,0)[bl]{\(A_{2}\)}}

\put(30,31){\oval(40,6)}


\put(00,19.5){\makebox(0,0)[bl]{\(\widetilde{B}_{3}:\)}}

\put(20,28){\vector(0,-1){4}}

\put(10.5,19.5){\makebox(0,0)[bl]{\(A_{1},A_{3},A_{7}\)}}

\put(30,21){\oval(40,6)}


\put(39.5,15){\oval(21,09)}

\put(30.5,13.5){\makebox(0,0)[bl]{\(\widetilde{B}_{3}\&\widetilde{B}_{4}:
A_{9}\)}}



\put(00,08.5){\makebox(0,0)[bl]{\(\widetilde{B}_{4}:\)}}

\put(20,18){\vector(0,-1){5}}

\put(13,08.5){\makebox(0,0)[bl]{\(A_{5},A_{8}\)}}

\put(30,10){\oval(40,6)}

\end{picture}
%
\begin{picture}(70,55)

\put(12,50){\makebox(0,0)[bl]{Table 2. Resultant priorities}}

\put(00,00){\line(1,0){70}} \put(00,38){\line(1,0){70}}
\put(00,48){\line(1,0){70}}

\put(00,0){\line(0,1){48}} \put(10,0){\line(0,1){48}}
\put(30,0){\line(0,1){48}} \put(50,0){\line(0,1){48}}
\put(70,0){\line(0,1){48}}

\put(01,34){\makebox(0,0)[bl]{\(A_{1}\)}}
\put(01,30){\makebox(0,0)[bl]{\(A_{2}\)}}
\put(01,26){\makebox(0,0)[bl]{\(A_{3}\)}}
\put(01,22){\makebox(0,0)[bl]{\(A_{4}\)}}
\put(01,18){\makebox(0,0)[bl]{\(A_{5}\)}}

\put(01,14){\makebox(0,0)[bl]{\(A_{6}\)}}
\put(01,10){\makebox(0,0)[bl]{\(A_{7}\)}}
\put(01,06){\makebox(0,0)[bl]{\(A_{8}\)}}
\put(01,02){\makebox(0,0)[bl]{\(A_{9}\)}}

\put(0.5,44){\makebox(0,0)[bl]{Alter-}}
\put(0.5,40){\makebox(0,0)[bl]{native}}

\put(13,44){\makebox(0,0)[bl]{Linear}}
\put(13,40){\makebox(0,0)[bl]{ordering}}

\put(33.5,44){\makebox(0,0)[bl]{Ranking}}

\put(53.5,44){\makebox(0,0)[bl]{``Fuzzy''}}
\put(53.5,40){\makebox(0,0)[bl]{ranking}}


\put(19,34){\makebox(0,0)[bl]{\(4\)}}
\put(39,34){\makebox(0,0)[bl]{\(3\)}}
\put(59,34){\makebox(0,0)[bl]{\(3\)}}

\put(19,30){\makebox(0,0)[bl]{\(3\)}}
\put(39,30){\makebox(0,0)[bl]{\(2\)}}
\put(59,30){\makebox(0,0)[bl]{\(2\)}}

\put(19,26){\makebox(0,0)[bl]{\(5\)}}
\put(39,26){\makebox(0,0)[bl]{\(3\)}}
\put(59,26){\makebox(0,0)[bl]{\(3\)}}

\put(19,22){\makebox(0,0)[bl]{\(1\)}}
\put(39,22){\makebox(0,0)[bl]{\(1\)}}
\put(59,22){\makebox(0,0)[bl]{\(1\)}}

\put(19,18){\makebox(0,0)[bl]{\(8\)}}
\put(39,18){\makebox(0,0)[bl]{\(4\)}}
\put(59,18){\makebox(0,0)[bl]{\(4\)}}

\put(19,14){\makebox(0,0)[bl]{\(2\)}}
\put(39,14){\makebox(0,0)[bl]{\(2\)}}
\put(56.5,13.4){\makebox(0,0)[bl]{\([1,2]\)}}

\put(19,10){\makebox(0,0)[bl]{\(6\)}}
\put(39,10){\makebox(0,0)[bl]{\(3\)}}
\put(59,10){\makebox(0,0)[bl]{\(3\)}}

\put(19,06){\makebox(0,0)[bl]{\(9\)}}
\put(39,06){\makebox(0,0)[bl]{\(4\)}}
\put(59,06){\makebox(0,0)[bl]{\(4\)}}

\put(19,02){\makebox(0,0)[bl]{\(7\)}}
\put(39,02){\makebox(0,0)[bl]{\(4\)}}
\put(56.5,01.4){\makebox(0,0)[bl]{\([3,4]\)}}

\end{picture}
\end{center}

\begin{center}
\begin{picture}(132,27)
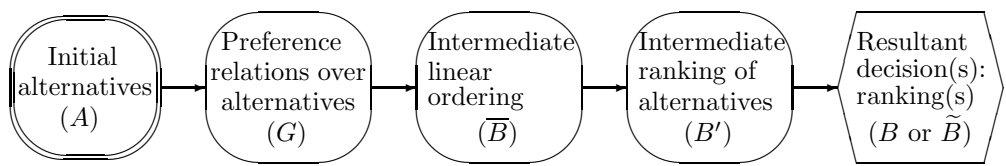
\put(14,00){\makebox(0,0)[bl] {Fig. 11. Decision making as
 preferences processing
 (\cite{lev86},\cite{lev98},\cite{levmih88})
}}


\put(05,19){\makebox(0,0)[bl]{Initial}}
\put(01,15){\makebox(0,0)[bl]{alternatives}}
\put(06.5,10){\makebox(0,0)[bl]{(\(A\))}}

\put(10,16){\oval(20,19)} \put(10,16){\oval(19,18)}

\put(20,16){\vector(1,0){6}}


\put(28,21){\makebox(0,0)[bl]{Preference}}
\put(26.6,17){\makebox(0,0)[bl]{relations over}}
\put(28,13){\makebox(0,0)[bl]{alternatives}}
\put(34,8){\makebox(0,0)[bl]{(\(G\))}}

\put(37,16){\oval(22,20)}

\put(48,16){\vector(1,0){6}}


\put(55.5,21){\makebox(0,0)[bl]{Intermediate}}
\put(55.5,17){\makebox(0,0)[bl]{linear }}
\put(55.5,13){\makebox(0,0)[bl]{ordering}}
\put(62,8){\makebox(0,0)[bl]{(\(\overline{B}\))}}

\put(65,16){\oval(22,20)}

\put(76,16){\vector(1,0){6}}


\put(83.5,21){\makebox(0,0)[bl]{Intermediate}}
\put(83.5,16.5){\makebox(0,0)[bl]{ranking of}}
\put(83.5,13){\makebox(0,0)[bl]{alternatives}}
\put(89.5,8){\makebox(0,0)[bl]{(\(B'\))}}

\put(93,16){\oval(22,20)}

\put(104,16){\vector(1,0){6}}


\put(112.5,21){\makebox(0,0)[bl]{Resultant}}
\put(112.5,16.5){\makebox(0,0)[bl]{decision(s):}}
\put(112.5,13){\makebox(0,0)[bl]{ranking(s)}}
\put(114,8){\makebox(0,0)[bl]{(\(B\) or \(\widetilde{B}\))}}

\put(110,16){\line(1,4){2.5}} \put(110,16){\line(1,-4){2.5}}
\put(132,16){\line(-1,4){2.5}} \put(132,16){\line(-1,-4){2.5}}
\put(112.5,26){\line(1,0){17}} \put(112.5,06){\line(1,0){17}}

\end{picture}
\end{center}

\subsection{Towards Reconfigurable Problem Solving Framework}

 Here a generalized problem solving framework is examined.
 An algorithm system (e.g., for data processing)
 is considered as a basic example.

 The basic parts of the framework are the following:
 (a) problem (problem structuring/formulation)
 (b) algorithm (operational part)
 (c) data part (obtaining data, elicitation of preferences, ...)
 (d) human part (preliminary learning, learning based on solving process)

 First, a statical structure of the algorithm system is analyzed.
 The following versions the algorithm/algorithm systems
 can be considered:

 1. usage of the basic algorithm(s) (Fig. 12);

 2. selection of the best algorithm from an algorithm base
 (while taking int account an input information)
 and its usage (Fig. 13);

 3. modification of the basic algorithm(s)
 (while taking int account an input information)
  and its usage (Fig. 14); and

  4. design of the new  algorithm(s)
 (while taking into account an input information)
  and its usage (Fig. 15).

 A simplified hierarchy of solving process components
 is depicted in Fig. 16.

 In general,
 the following problem solving frameworks can be considered:

 1. One-phase framework: problem solving.

 2. Two-phase framework:
 (i) problem structuring/formulation and
 (ii) problem solving.

 3. Three-phase framework:
  (i) problem structuring/formulation;
 (ii) problem solving; and
  (iii) analysis of results.

  4. Adaptive three-phase framework with feedback:
  (a) problem structuring/formulation;
 (b) problem solving (including analysis of intermediate results and
 and problem reformulation and resolving);
  (c) analysis of results
 (including analysis of results
 and problem reformulation and resolving).

\begin{center}
\begin{picture}(70,33)

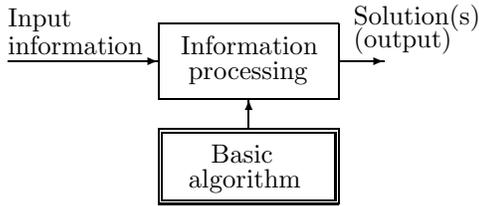
\put(03,00){\makebox(0,0)[bl]{Fig. 12. Usage of a basic
algorithm}}

\put(0,29){\makebox(0,0)[bl]{Input}}
\put(0,26){\makebox(0,0)[bl]{information}}
\put(0,25){\vector(1,0){20}}


\put(20,30){\line(1,0){24}} \put(20,20){\line(1,0){24}}
\put(20,20){\line(0,1){10}} \put(44,20){\line(0,1){10}}

\put(23,26){\makebox(0,0)[bl]{Information}}
\put(24,22){\makebox(0,0)[bl]{processing}}

\put(44,25){\vector(1,0){6}}

\put(46,29){\makebox(0,0)[bl]{Solution(s)}}
\put(46,26){\makebox(0,0)[bl]{(output)}}


\put(20,16){\line(1,0){24}} \put(20,06){\line(1,0){24}}
\put(20,06){\line(0,1){10}} \put(44,06){\line(0,1){10}}

\put(20.5,15.5){\line(1,0){23}} \put(20.5,6.5){\line(1,0){23}}
\put(20.5,6.5){\line(0,1){9}} \put(43.5,6.5){\line(0,1){9}}

\put(27,11.5){\makebox(0,0)[bl]{Basic }}
\put(24,07.5){\makebox(0,0)[bl]{algorithm}}

\put(32,16){\vector(0,1){4}}

\end{picture}
%
\begin{picture}(70,47)

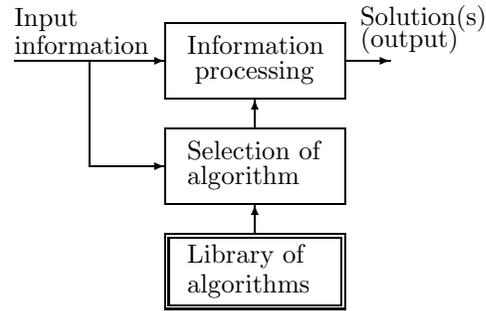
\put(01,00){\makebox(0,0)[bl]{Fig. 13. Usage of a selected
algorithm}}

\put(0,43){\makebox(0,0)[bl]{Input}}
\put(0,40){\makebox(0,0)[bl]{information}}
\put(0,39){\vector(1,0){20}}

\put(10,39){\line(0,-1){14}}\put(10,25){\vector(1,0){10}}


\put(20,44){\line(1,0){24}} \put(20,34){\line(1,0){24}}
\put(20,34){\line(0,1){10}} \put(44,34){\line(0,1){10}}

\put(23,40){\makebox(0,0)[bl]{Information}}
\put(24,36){\makebox(0,0)[bl]{processing}}

\put(44,39){\vector(1,0){6}}

\put(46,43){\makebox(0,0)[bl]{Solution(s)}}
\put(46,40){\makebox(0,0)[bl]{(output)}}


\put(20,30){\line(1,0){24}} \put(20,20){\line(1,0){24}}
\put(20,20){\line(0,1){10}} \put(44,20){\line(0,1){10}}

\put(23,26){\makebox(0,0)[bl]{Selection of }}
\put(23,22){\makebox(0,0)[bl]{algorithm}}

\put(32,30){\vector(0,1){4}}


\put(20,16){\line(1,0){24}} \put(20,06){\line(1,0){24}}
\put(20,06){\line(0,1){10}} \put(44,06){\line(0,1){10}}

\put(20.5,15.5){\line(1,0){23}} \put(20.5,6.5){\line(1,0){23}}
\put(20.5,6.5){\line(0,1){9}} \put(43.5,6.5){\line(0,1){9}}

\put(23,11.5){\makebox(0,0)[bl]{Library of }}
\put(23,07.5){\makebox(0,0)[bl]{algorithms}}

\put(32,16){\vector(0,1){4}}

\end{picture}
\end{center}

\begin{center}
\begin{picture}(70,47)

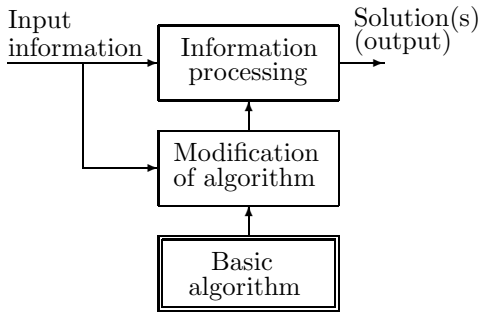
\put(00,00){\makebox(0,0)[bl]{Fig. 14. Modification a basic
algorithm}}

\put(0,43){\makebox(0,0)[bl]{Input}}
\put(0,40){\makebox(0,0)[bl]{information}}
\put(0,39){\vector(1,0){20}}

\put(10,39){\line(0,-1){14}}\put(10,25){\vector(1,0){10}}


\put(20,44){\line(1,0){24}} \put(20,34){\line(1,0){24}}
\put(20,34){\line(0,1){10}} \put(44,34){\line(0,1){10}}

\put(23,40){\makebox(0,0)[bl]{Information}}
\put(24,36){\makebox(0,0)[bl]{processing}}

\put(44,39){\vector(1,0){6}}

\put(46,43){\makebox(0,0)[bl]{Solution(s)}}
\put(46,40){\makebox(0,0)[bl]{(output)}}


\put(20,30){\line(1,0){24}} \put(20,20){\line(1,0){24}}
\put(20,20){\line(0,1){10}} \put(44,20){\line(0,1){10}}

\put(22,26){\makebox(0,0)[bl]{Modification }}
\put(22,22){\makebox(0,0)[bl]{of algorithm}}

\put(32,30){\vector(0,1){4}}


\put(20,16){\line(1,0){24}} \put(20,06){\line(1,0){24}}
\put(20,06){\line(0,1){10}} \put(44,06){\line(0,1){10}}

\put(20.5,15.5){\line(1,0){23}} \put(20.5,6.5){\line(1,0){23}}
\put(20.5,6.5){\line(0,1){9}} \put(43.5,6.5){\line(0,1){9}}

\put(27,11.5){\makebox(0,0)[bl]{Basic }}
\put(24,07.5){\makebox(0,0)[bl]{algorithm}}

\put(32,16){\vector(0,1){4}}

\end{picture}
%
\begin{picture}(70,47)

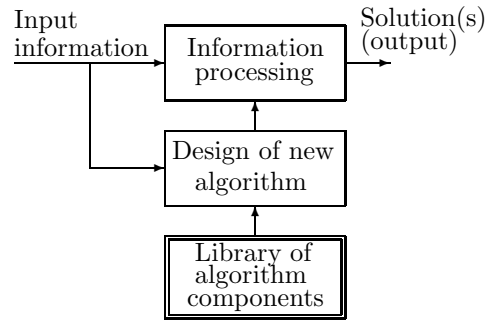
\put(002,00){\makebox(0,0)[bl]{Fig. 15. Design of a new
algorithm}}

\put(0,43){\makebox(0,0)[bl]{Input}}
\put(0,40){\makebox(0,0)[bl]{information}}
\put(0,39){\vector(1,0){20}}

\put(10,39){\line(0,-1){14}}\put(10,25){\vector(1,0){10}}


\put(20,44){\line(1,0){24}} \put(20,34){\line(1,0){24}}
\put(20,34){\line(0,1){10}} \put(44,34){\line(0,1){10}}

\put(23,40){\makebox(0,0)[bl]{Information}}
\put(24,36){\makebox(0,0)[bl]{processing}}

\put(44,39){\vector(1,0){6}}

\put(46,43){\makebox(0,0)[bl]{Solution(s)}}
\put(46,40){\makebox(0,0)[bl]{(output)}}


\put(20,30){\line(1,0){24}} \put(20,20){\line(1,0){24}}
\put(20,20){\line(0,1){10}} \put(44,20){\line(0,1){10}}

\put(21,25.5){\makebox(0,0)[bl]{Design of new  }}
\put(24,21.5){\makebox(0,0)[bl]{algorithm}}

\put(32,30){\vector(0,1){4}}


\put(20,16){\line(1,0){24}} \put(20,05){\line(1,0){24}}
\put(20,05){\line(0,1){11}} \put(44,05){\line(0,1){11}}

\put(20.5,15.5){\line(1,0){23}} \put(20.5,5.5){\line(1,0){23}}
\put(20.5,5.5){\line(0,1){10}} \put(43.5,5.5){\line(0,1){10}}

\put(24,12){\makebox(0,0)[bl]{Library of}}
\put(24,09){\makebox(0,0)[bl]{algorithm }}
\put(23,06){\makebox(0,0)[bl]{components}}

\put(32,16){\vector(0,1){4}}

\end{picture}
\end{center}

\begin{center}
\begin{picture}(98,64)
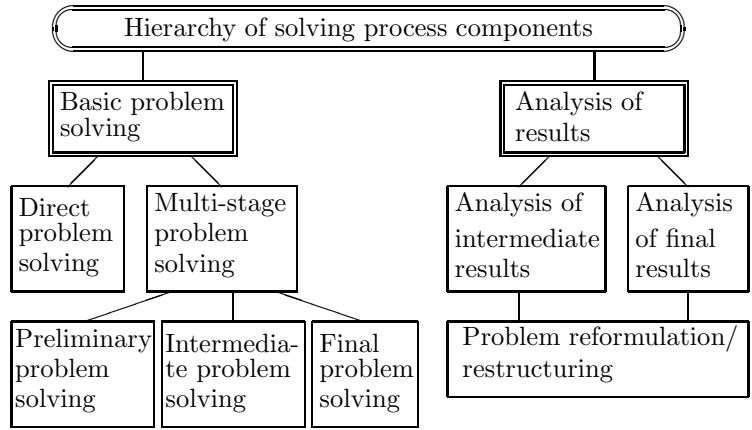
\put(07,00){\makebox(0,0)[bl]{Fig. 16. Problem solving with
 analysis, reformulation}}

\put(47.5,59){\oval(84,6)} \put(47.5,59){\oval(83,5)}

\put(15,57.3){\makebox(0,0)[bl]{Hierarchy of solving process
components}}

\put(17.5,56){\line(0,-1){4}} \put(77.5,56){\line(0,-1){4}}


\put(5,42){\line(1,0){25}} \put(5,52){\line(1,0){25}}
\put(5,42){\line(0,1){10}} \put(30,42){\line(0,1){10}}

\put(5.5,42.5){\line(1,0){24}} \put(5.5,51.5){\line(1,0){24}}
\put(5.5,42.5){\line(0,1){09}} \put(29.5,42.5){\line(0,1){09}}

\put(06.5,47.5){\makebox(0,0)[bl]{Basic problem}}
\put(06.5,44){\makebox(0,0)[bl]{solving}}

\put(7.5,38){\line(1,1){4}} \put(28,38){\line(-1,1){4}}


\put(00,24){\line(1,0){15}} \put(00,38){\line(1,0){15}}
\put(00,24){\line(0,1){14}} \put(15,24){\line(0,1){14}}

\put(01,34){\makebox(0,0)[bl]{Direct }}
\put(01,30){\makebox(0,0)[bl]{problem}}
\put(01,26){\makebox(0,0)[bl]{solving}}


\put(18,24){\line(1,0){20}} \put(18,38){\line(1,0){20}}
\put(18,24){\line(0,1){14}} \put(38,24){\line(0,1){14}}

\put(19,34){\makebox(0,0)[bl]{Multi-stage}}
\put(19,30){\makebox(0,0)[bl]{problem}}
\put(19,26){\makebox(0,0)[bl]{solving}}


\put(00,06){\line(1,0){19}} \put(00,20){\line(1,0){19}}
\put(00,06){\line(0,1){14}} \put(19,06){\line(0,1){14}}

\put(0.5,16){\makebox(0,0)[bl]{Preliminary}}
\put(0.5,12){\makebox(0,0)[bl]{problem}}
\put(0.5,08){\makebox(0,0)[bl]{solving}}

\put(9.5,20){\line(3,1){12}}


\put(20,06){\line(1,0){19}} \put(20,20){\line(1,0){19}}
\put(20,06){\line(0,1){14}} \put(39,06){\line(0,1){14}}

\put(20.5,16){\makebox(0,0)[bl]{Intermedia-}}
\put(20.5,12){\makebox(0,0)[bl]{te problem}}
\put(20.5,08){\makebox(0,0)[bl]{solving}}

\put(29.5,20){\line(0,1){4}}


\put(40,06){\line(1,0){14}} \put(40,20){\line(1,0){14}}
\put(40,06){\line(0,1){14}} \put(54,06){\line(0,1){14}}

\put(41,16){\makebox(0,0)[bl]{Final}}
\put(41,12){\makebox(0,0)[bl]{problem}}
\put(41,08){\makebox(0,0)[bl]{solving}}

\put(47,20){\line(-3,1){12}}


\put(65,42){\line(1,0){25}} \put(65,52){\line(1,0){25}}
\put(65,42){\line(0,1){10}} \put(90,42){\line(0,1){10}}

\put(65.5,42.5){\line(1,0){24}} \put(65.5,51.5){\line(1,0){24}}
\put(65.5,42.5){\line(0,1){09}} \put(89.5,42.5){\line(0,1){09}}

\put(67,47.5){\makebox(0,0)[bl]{Analysis of}}
\put(67,44){\makebox(0,0)[bl]{results}}

\put(67.5,38){\line(1,1){4}} \put(90,38){\line(-1,1){4}}

\put(67.5,24){\line(0,-1){4}} \put(90,24){\line(0,-1){4}}


\put(58,24){\line(1,0){21}} \put(58,38){\line(1,0){21}}
\put(58,24){\line(0,1){14}} \put(79,24){\line(0,1){14}}

\put(59,34){\makebox(0,0)[bl]{Analysis of}}
\put(59,30){\makebox(0,0)[bl]{intermediate}}
\put(59,26){\makebox(0,0)[bl]{results}}


\put(82,24){\line(1,0){16}} \put(82,38){\line(1,0){16}}
\put(82,24){\line(0,1){14}} \put(98,24){\line(0,1){14}}

\put(83,34){\makebox(0,0)[bl]{Analysis}}
\put(83,30){\makebox(0,0)[bl]{of final}}
\put(83,26){\makebox(0,0)[bl]{results}}


\put(58,10){\line(1,0){40}} \put(58,20){\line(1,0){40}}
\put(58,10){\line(0,1){10}} \put(98,10){\line(0,1){10}}

\put(60,16){\makebox(0,0)[bl]{Problem reformulation/}}
\put(60,12){\makebox(0,0)[bl]{restructuring}}

\end{picture}
\end{center}

 The following three illustrative examples of multistage solving strategies
 are presented for set consisting of nine initial alternatives (Table 1):
 (i) two-stage series ranking strategy with preliminary linear ordering (Fig. 17);
 (ii) two-stage series ranking strategy with preliminary preferences (Fig. 18);
 (iii) three-stage series ranking strategy (Fig. 19); and
 (iv) three-stage series-parallel strategy (with aggregation) (Fig. 20).

 Fig. 21 depicts  adaptive three-stage framework for problem solving.
 In addition, it is possible to consider a solving framework in
 which problem formulation/structuring is executed during the
 solving process for a preliminary problem.


 The stage
 ``Reconfiguration of problem solving process''
 can be considered as the following::
 (i) selection of another algorithm (solving scheme) (Fig. 13),
 (ii) modification of the algorithm (solving scheme) (Fig. 14),
 (iii) design of a new algorithm (solving scheme) (Fig. 15).
 Concurrently, information operations
 (searching/acquisition/usage)
 can be used
 (e.g., new/additional preferences,
 new/additional reference points of intermediate decisions).

 In the case of a man-machine procedure,
 it is possible to examine change or re-organization
 of an expert team and types (e.g., mode, support procedure(s)) of man-machine
 interaction.

 The stage
 ``Problem reformulation/restructuring''
 can be considered as the following::
 (i) change of the initial problem,
 (ii) modification of the initial problem
 (e.g., change of the problem parameters, initial data),
 (iii) building of a new problem framework (i.e., a composite
 problem).

 Fig. 22 depicts an example to  illustrate  a series process with problem
 reformulation.

\begin{center}
\begin{picture}(58,67)

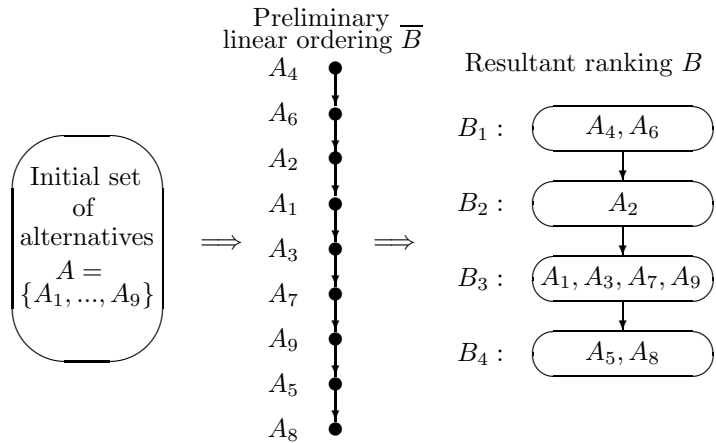
\put(00,00){\makebox(0,0)[bl]{Fig. 17. Two-stage series ranking
strategy (with linear ordering)}}

\put(02.4,42.5){\makebox(0,0)[bl]{Initial set}}
\put(01,38.3){\makebox(0,0)[bl]{~~~~~ of}}
\put(01,34.5){\makebox(0,0)[bl]{alternatives}}

\put(01,30){\makebox(0,0)[bl]{~~~ \(A = \)}}

\put(01,26){\makebox(0,0)[bl]{\(\{A_{1},...,A_{9}\}\)}}

\put(10,34){\oval(20,30)}


\put(25,34){\makebox(0,0)[bl]{\(\Longrightarrow \)}}


\put(32,63){\makebox(0,0)[bl]{Preliminary}}
\put(28,60){\makebox(0,0)[bl]{linear ordering \(\overline{B}\)}}


\put(34,56.4){\makebox(0,0)[bl]{\(A_{4}\)}}

\put(43,58){\circle*{1.8}}


\put(43,58){\vector(0,-1){5}}

\put(34,50.4){\makebox(0,0)[bl]{\(A_{6}\)}}

\put(43,52){\circle*{1.8}}


\put(43,52){\vector(0,-1){5}}

\put(34,44.4){\makebox(0,0)[bl]{\(A_{2}\)}}

\put(43,46){\circle*{1.8}}


\put(43,46){\vector(0,-1){5}}

\put(34,38.4){\makebox(0,0)[bl]{\(A_{1}\)}}

\put(43,40){\circle*{1.8}}


\put(43,40){\vector(0,-1){5}}

\put(34,32.4){\makebox(0,0)[bl]{\(A_{3}\)}}

\put(43,34){\circle*{1.8}}


\put(43,34){\vector(0,-1){5}}

\put(34,26.4){\makebox(0,0)[bl]{\(A_{7}\)}}

\put(43,28){\circle*{1.8}}


\put(43,28){\vector(0,-1){5}}

\put(34,20.4){\makebox(0,0)[bl]{\(A_{9}\)}}

\put(43,22){\circle*{1.8}}


\put(43,22){\vector(0,-1){5}}

\put(34,14.4){\makebox(0,0)[bl]{\(A_{5}\)}}

\put(43,16){\circle*{1.8}}


\put(43,16){\vector(0,-1){5}}

\put(34,08.4){\makebox(0,0)[bl]{\(A_{8}\)}}

\put(43,10){\circle*{1.8}}



\put(48,34){\makebox(0,0)[bl]{\(\Longrightarrow \)}}


\end{picture}
%
\begin{picture}(42,50)


\put(01,57){\makebox(0,0)[bl]{Resultant ranking \(B\)}}


\put(00,48.5){\makebox(0,0)[bl]{\(B_{1}:\)}}

\put(17,48.5){\makebox(0,0)[bl]{\(A_{4},A_{6}\)}}

\put(22,50){\oval(24,6)}


\put(00,38.5){\makebox(0,0)[bl]{\(B_{2}:\)}}

\put(22,47){\vector(0,-1){4}}

\put(20,38.5){\makebox(0,0)[bl]{\(A_{2}\)}}

\put(22,40){\oval(24,6)}


\put(00,28.5){\makebox(0,0)[bl]{\(B_{3}:\)}}

\put(22,37){\vector(0,-1){4}}

\put(11,28.5){\makebox(0,0)[bl]{\(A_{1},A_{3},A_{7},A_{9}\)}}

\put(22,30){\oval(24,6)}


\put(00,18.5){\makebox(0,0)[bl]{\(B_{4}:\)}}

\put(22,27){\vector(0,-1){4}}

\put(17,18.5){\makebox(0,0)[bl]{\(A_{5},A_{8}\)}}

\put(22,20){\oval(24,6)}

\end{picture}
\end{center}

\begin{center}
\begin{picture}(58,52)

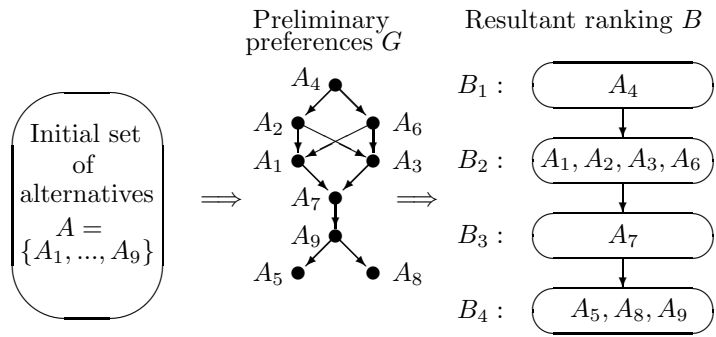
\put(00,00){\makebox(0,0)[bl]{Fig. 18. Two-stage series ranking
strategy (with preferences)}}

\put(02.4,32.5){\makebox(0,0)[bl]{Initial set}}
\put(01,28.3){\makebox(0,0)[bl]{~~~~~ of}}
\put(01,24.5){\makebox(0,0)[bl]{alternatives}}

\put(01,20){\makebox(0,0)[bl]{~~~ \(A = \)}}

\put(01,16){\makebox(0,0)[bl]{\(\{A_{1},...,A_{9}\}\)}}

\put(10,24){\oval(20,30)}


\put(25,24){\makebox(0,0)[bl]{\(\Longrightarrow \)}}


\put(32,47){\makebox(0,0)[bl]{Preliminary}}
\put(31,44){\makebox(0,0)[bl]{preferences \(G\)}}


\put(37,39.4){\makebox(0,0)[bl]{\(A_{4}\)}}

\put(43,40){\circle*{1.8}}

\put(43,40){\vector(-1,-1){4}} \put(43,40){\vector(1,-1){4}}

\put(32,33.4){\makebox(0,0)[bl]{\(A_{2}\)}}

\put(38,35){\circle*{1.8}}

\put(50.5,33.4){\makebox(0,0)[bl]{\(A_{6}\)}}

\put(48,35){\circle*{1.8}}

\put(38,35){\vector(0,-1){4}} \put(48,35){\vector(0,-1){4}}

\put(38,35){\vector(2,-1){9}} \put(48,35){\vector(-2,-1){9}}

\put(32,28.4){\makebox(0,0)[bl]{\(A_{1}\)}}

\put(38,30){\circle*{1.8}}

\put(50.5,28.4){\makebox(0,0)[bl]{\(A_{3}\)}}

\put(48,30){\circle*{1.8}}

\put(38,30){\vector(1,-1){4}} \put(48,30){\vector(-1,-1){4}}


\put(37,23.4){\makebox(0,0)[bl]{\(A_{7}\)}}

\put(43,25){\circle*{1.8}}

\put(43,25){\vector(0,-1){4}}


\put(37,18.4){\makebox(0,0)[bl]{\(A_{9}\)}}

\put(43,20){\circle*{1.8}}

\put(43,20){\vector(-1,-1){4}} \put(43,20){\vector(1,-1){4}}

\put(32,13.4){\makebox(0,0)[bl]{\(A_{5}\)}}

\put(38,15){\circle*{1.8}}

\put(50.5,13.4){\makebox(0,0)[bl]{\(A_{8}\)}}

\put(48,15){\circle*{1.8}}


\put(51,24){\makebox(0,0)[bl]{\(\Longrightarrow \)}}


\end{picture}
\begin{picture}(42,50)


\put(01,47){\makebox(0,0)[bl]{Resultant ranking \(B\)}}


\put(00,38.5){\makebox(0,0)[bl]{\(B_{1}:\)}}

\put(20,38.5){\makebox(0,0)[bl]{\(A_{4}\)}}

\put(22,40){\oval(24,6)}


\put(00,28.5){\makebox(0,0)[bl]{\(B_{2}:\)}}

\put(22,37){\vector(0,-1){4}}

\put(11,28.5){\makebox(0,0)[bl]{\(A_{1},A_{2},A_{3},A_{6}\)}}

\put(22,30){\oval(24,6)}


\put(00,18.5){\makebox(0,0)[bl]{\(B_{3}:\)}}

\put(22,27){\vector(0,-1){4}}

\put(20,18.5){\makebox(0,0)[bl]{\(A_{7}\)}}

\put(22,20){\oval(24,6)}


\put(00,08.5){\makebox(0,0)[bl]{\(B_{4}:\)}}

\put(22,17){\vector(0,-1){4}}

\put(15,8.5){\makebox(0,0)[bl]{\(A_{5},A_{8},A_{9}\)}}

\put(22,10){\oval(24,6)}

\end{picture}
\end{center}

\begin{center}
\begin{picture}(50,71)

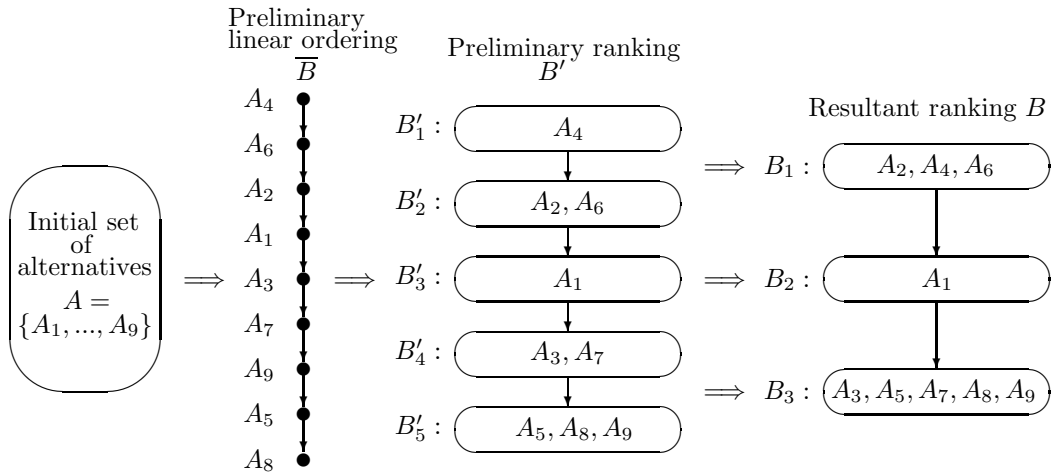
\put(28,00){\makebox(0,0)[bl]{Fig. 19. Three-stage series ranking
strategy}}

\put(02.4,40.5){\makebox(0,0)[bl]{Initial set}}
\put(01,37.5){\makebox(0,0)[bl]{~~~~~ of}}

\put(01,34.5){\makebox(0,0)[bl]{alternatives}}

\put(07,30){\makebox(0,0)[bl]{\(A =\)}}

\put(01,26){\makebox(0,0)[bl]{\(\{A_{1},...,A_{9}\}\)}}

\put(10,34){\oval(20,30)}


\put(23,33){\makebox(0,0)[bl]{\(\Longrightarrow \)}}


\put(29,67){\makebox(0,0)[bl]{Preliminary}}
\put(29,64){\makebox(0,0)[bl]{linear ordering}}
\put(38,60.5){\makebox(0,0)[bl]{\(\overline{B}\)}}


\put(31,56.4){\makebox(0,0)[bl]{\(A_{4}\)}}

\put(39,58){\circle*{1.8}}


\put(39,58){\vector(0,-1){5}}

\put(31,50.4){\makebox(0,0)[bl]{\(A_{6}\)}}

\put(39,52){\circle*{1.8}}


\put(39,52){\vector(0,-1){5}}

\put(31,44.4){\makebox(0,0)[bl]{\(A_{2}\)}}

\put(39,46){\circle*{1.8}}


\put(39,46){\vector(0,-1){5}}

\put(31,38.4){\makebox(0,0)[bl]{\(A_{1}\)}}

\put(39,40){\circle*{1.8}}


\put(39,40){\vector(0,-1){5}}

\put(31,32.4){\makebox(0,0)[bl]{\(A_{3}\)}}

\put(39,34){\circle*{1.8}}


\put(39,34){\vector(0,-1){5}}

\put(31,26.4){\makebox(0,0)[bl]{\(A_{7}\)}}

\put(39,28){\circle*{1.8}}


\put(39,28){\vector(0,-1){5}}

\put(31,20.4){\makebox(0,0)[bl]{\(A_{9}\)}}

\put(39,22){\circle*{1.8}}


\put(39,22){\vector(0,-1){5}}

\put(31,14.4){\makebox(0,0)[bl]{\(A_{5}\)}}

\put(39,16){\circle*{1.8}}


\put(39,16){\vector(0,-1){5}}

\put(31,08.4){\makebox(0,0)[bl]{\(A_{8}\)}}

\put(39,10){\circle*{1.8}}



\put(43,33){\makebox(0,0)[bl]{\(\Longrightarrow \)}}


\end{picture}
%
\begin{picture}(87,60)


\put(07,63){\makebox(0,0)[bl]{Preliminary ranking}}
\put(19,60.5){\makebox(0,0)[bl]{\(B'\)}}


\put(00,52.5){\makebox(0,0)[bl]{\(B'_{1}:\)}}

\put(21,52.5){\makebox(0,0)[bl]{\(A_{4}\)}}

\put(23,54){\oval(30,6)}


\put(00,42.5){\makebox(0,0)[bl]{\(B'_{2}:\)}}

\put(23,51){\vector(0,-1){4}}

\put(18,42.5){\makebox(0,0)[bl]{\(A_{2},A_{6}\)}}

\put(23,44){\oval(30,6)}


\put(00,32.5){\makebox(0,0)[bl]{\(B'_{3}:\)}}

\put(23,41){\vector(0,-1){4}}

\put(21,32.5){\makebox(0,0)[bl]{\(A_{1}\)}}

\put(23,34){\oval(30,6)}


\put(00,22.5){\makebox(0,0)[bl]{\(B'_{4}:\)}}

\put(23,31){\vector(0,-1){4}}

\put(18,22.5){\makebox(0,0)[bl]{\(A_{3},A_{7}\)}}

\put(23,24){\oval(30,6)}


\put(00,12.5){\makebox(0,0)[bl]{\(B'_{5}:\)}}

\put(23,21){\vector(0,-1){4}}

\put(16,12.5){\makebox(0,0)[bl]{\(A_{5},A_{8},A_{9}\)}}

\put(23,14){\oval(30,6)}


\put(41,48){\makebox(0,0)[bl]{\(\Longrightarrow \)}}
\put(41,33){\makebox(0,0)[bl]{\(\Longrightarrow \)}}
\put(41,18){\makebox(0,0)[bl]{\(\Longrightarrow \)}}


\put(55,55){\makebox(0,0)[bl]{Resultant ranking \(B\)}}


\put(49,47.5){\makebox(0,0)[bl]{\(B_{1}:\)}}


\put(64,47.5){\makebox(0,0)[bl]{\(A_{2},A_{4},A_{6}\)}}

\put(72,49){\oval(30,6)}


\put(49,32.5){\makebox(0,0)[bl]{\(B_{2}:\)}}

\put(72,46){\vector(0,-1){9}}

\put(70,32.5){\makebox(0,0)[bl]{\(A_{1}\)}}

\put(72,34){\oval(30,6)}


\put(49,17.5){\makebox(0,0)[bl]{\(B_{3}:\)}}

\put(72,31){\vector(0,-1){9}}

\put(58,17.5){\makebox(0,0)[bl]{\(A_{3},A_{5},A_{7},A_{8},A_{9}\)}}

\put(72,19){\oval(30,6)}


\end{picture}
\end{center}
\begin{center}
\begin{picture}(50,103)

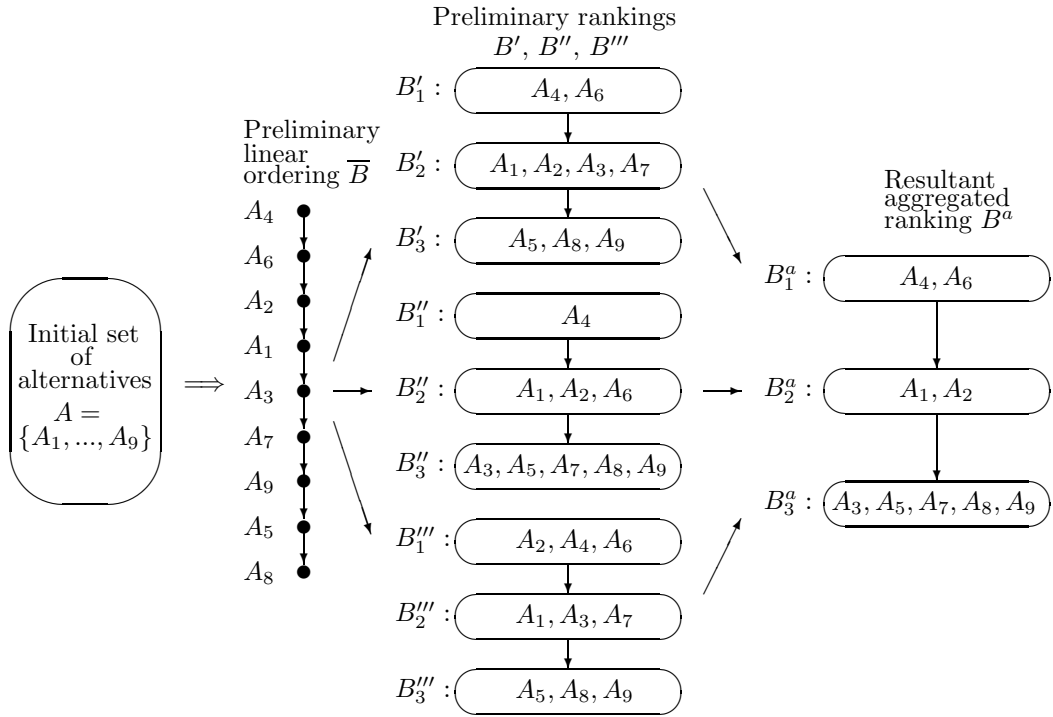
\put(13,00){\makebox(0,0)[bl]{Fig. 20. Three-stage
 series-parallel ranking strategy (with aggregation)}}

\put(02.4,56.5){\makebox(0,0)[bl]{Initial set}}
\put(01,53.5){\makebox(0,0)[bl]{~~~~~ of}}

\put(01,50.5){\makebox(0,0)[bl]{alternatives}}

\put(01,46){\makebox(0,0)[bl]{~~~ \(A =\)}}

\put(01,42){\makebox(0,0)[bl]{\(\{A_{1},...,A_{9}\}\)}}

\put(10,50){\oval(20,30)}


\put(23,50){\makebox(0,0)[bl]{\(\Longrightarrow \)}}


\put(31,83){\makebox(0,0)[bl]{Preliminary}}
\put(31,80.5){\makebox(0,0)[bl]{linear}}
\put(31,77){\makebox(0,0)[bl]{ordering \(\overline{B}\)}}


\put(31,72.4){\makebox(0,0)[bl]{\(A_{4}\)}}

\put(39,74){\circle*{1.8}}


\put(39,74){\vector(0,-1){5}}

\put(31,66.4){\makebox(0,0)[bl]{\(A_{6}\)}}

\put(39,68){\circle*{1.8}}


\put(39,68){\vector(0,-1){5}}

\put(31,60.4){\makebox(0,0)[bl]{\(A_{2}\)}}

\put(39,62){\circle*{1.8}}


\put(39,62){\vector(0,-1){5}}

\put(31,54.4){\makebox(0,0)[bl]{\(A_{1}\)}}

\put(39,56){\circle*{1.8}}


\put(39,56){\vector(0,-1){5}}

\put(31,48.4){\makebox(0,0)[bl]{\(A_{3}\)}}

\put(39,50){\circle*{1.8}}


\put(39,50){\vector(0,-1){5}}

\put(31,42.4){\makebox(0,0)[bl]{\(A_{7}\)}}

\put(39,44){\circle*{1.8}}


\put(39,44){\vector(0,-1){5}}

\put(31,36.4){\makebox(0,0)[bl]{\(A_{9}\)}}

\put(39,38){\circle*{1.8}}


\put(39,38){\vector(0,-1){5}}

\put(31,30.4){\makebox(0,0)[bl]{\(A_{5}\)}}

\put(39,32){\circle*{1.8}}


\put(39,32){\vector(0,-1){5}}

\put(31,24.4){\makebox(0,0)[bl]{\(A_{8}\)}}

\put(39,26){\circle*{1.8}}



\put(43,54){\vector(1,3){5}}

\put(43,50){\vector(1,0){5}}

\put(43,46){\vector(1,-3){5}}



\end{picture}
%
\begin{picture}(87,103)


\put(05,98){\makebox(0,0)[bl]{Preliminary rankings}}
\put(13,94){\makebox(0,0)[bl]{\(B'\), \(B''\), \(B'''\)}}


\put(00,88.5){\makebox(0,0)[bl]{\(B'_{1}:\)}}

\put(18,88.5){\makebox(0,0)[bl]{\(A_{4},A_{6}\)}}

\put(23,90){\oval(30,6)}


\put(00,78.5){\makebox(0,0)[bl]{\(B'_{2}:\)}}

\put(23,87){\vector(0,-1){4}}

\put(12.5,78.5){\makebox(0,0)[bl]{\(A_{1},A_{2},A_{3},A_{7}\)}}

\put(23,80){\oval(30,6)}


\put(00,68.5){\makebox(0,0)[bl]{\(B'_{3}:\)}}

\put(23,77){\vector(0,-1){4}}

\put(15,68.5){\makebox(0,0)[bl]{\(A_{5},A_{8},A_{9}\)}}

\put(23,70){\oval(30,6)}


\put(00,58.5){\makebox(0,0)[bl]{\(B''_{1}:\)}}

\put(22,58.5){\makebox(0,0)[bl]{\(A_{4}\)}}

\put(23,60){\oval(30,6)}


\put(00,48.5){\makebox(0,0)[bl]{\(B''_{2}:\)}}

\put(23,57){\vector(0,-1){4}}

\put(16,48.5){\makebox(0,0)[bl]{\(A_{1},A_{2},A_{6}\)}}

\put(23,50){\oval(30,6)}


\put(00,38.5){\makebox(0,0)[bl]{\(B''_{3}:\)}}

\put(23,47){\vector(0,-1){4}}

\put(09,38.5){\makebox(0,0)[bl]{\(A_{3},A_{5},A_{7},A_{8},A_{9}\)}}

\put(23,40){\oval(30,6)}


\put(00,28.5){\makebox(0,0)[bl]{\(B'''_{1}:\)}}

\put(16,28.5){\makebox(0,0)[bl]{\(A_{2},A_{4},A_{6}\)}}

\put(23,30){\oval(30,6)}


\put(00,18.5){\makebox(0,0)[bl]{\(B'''_{2}:\)}}

\put(23,27){\vector(0,-1){4}}

\put(16,18.5){\makebox(0,0)[bl]{\(A_{1},A_{3},A_{7}\)}}

\put(23,20){\oval(30,6)}


\put(00,08.5){\makebox(0,0)[bl]{\(B'''_{3}:\)}}

\put(23,17){\vector(0,-1){4}}

\put(16,08.5){\makebox(0,0)[bl]{\(A_{5},A_{8},A_{9}\)}}

\put(23,10){\oval(30,6)}


\put(41,77){\vector(1,-2){5}}

\put(41,50){\vector(1,0){5}}

\put(41,23){\vector(1,2){5}}


\put(65,77){\makebox(0,0)[bl]{Resultant}}
\put(65,74){\makebox(0,0)[bl]{aggregated}}
\put(65,71){\makebox(0,0)[bl]{ranking \(B^{a}\)}}


\put(49,63.5){\makebox(0,0)[bl]{\(B^{a}_{1}:\)}}

\put(67,63.5){\makebox(0,0)[bl]{\(A_{4},A_{6}\)}}

\put(72,65){\oval(30,6)}


\put(49,48.5){\makebox(0,0)[bl]{\(B^{a}_{2}:\)}}

\put(72,62){\vector(0,-1){9}}

\put(67,48.5){\makebox(0,0)[bl]{\(A_{1},A_{2}\)}}

\put(72,50){\oval(30,6)}


\put(49,33.5){\makebox(0,0)[bl]{\(B^{a}_{3}:\)}}

\put(72,47){\vector(0,-1){9}}

\put(58,33.5){\makebox(0,0)[bl]{\(A_{3},A_{5},A_{7},A_{8},A_{9}\)}}

\put(72,35){\oval(30,6)}


\end{picture}
\end{center}

\begin{center}
\begin{picture}(81,70)
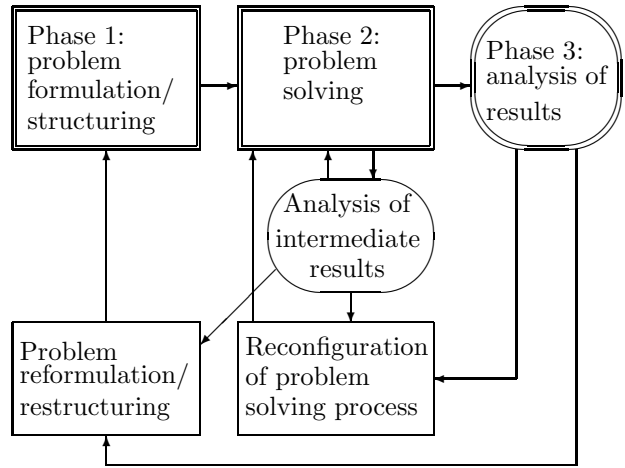
\put(02.5,00){\makebox(0,0)[bl]{Fig. 21. Problem solving with
 reconfiguration}}


\put(00,48){\line(1,0){25}} \put(00,67){\line(1,0){25}}
\put(00,48){\line(0,1){19}} \put(25,48){\line(0,1){19}}

\put(0.5,48.5){\line(1,0){24}} \put(0.5,66.5){\line(1,0){24}}
\put(0.5,48.5){\line(0,1){18}} \put(24.5,48.5){\line(0,1){18}}

\put(25,56.5){\vector(1,0){5}}

\put(02,62){\makebox(0,0)[bl]{Phase 1:}}
\put(02,58){\makebox(0,0)[bl]{problem}}
\put(02,54){\makebox(0,0)[bl]{formulation/}}
\put(02,50){\makebox(0,0)[bl]{structuring}}


\put(30,48){\line(1,0){26}} \put(30,67){\line(1,0){26}}
\put(30,48){\line(0,1){19}} \put(56,48){\line(0,1){19}}

\put(30.5,48.5){\line(1,0){25}} \put(30.5,66.5){\line(1,0){25}}
\put(30.5,48.5){\line(0,1){18}} \put(55.5,48.5){\line(0,1){18}}

\put(56,56.5){\vector(1,0){5}}

\put(36,62){\makebox(0,0)[bl]{Phase 2:}}
\put(36,58){\makebox(0,0)[bl]{problem}}
\put(36,54){\makebox(0,0)[bl]{solving}}

\put(48,48){\vector(0,-1){4}}

\put(42,44){\vector(0,1){4}}


\put(71,57.5){\oval(20,19)} \put(71,57.5){\oval(19,18)}

\put(63,60){\makebox(0,0)[bl]{Phase 3:}}
\put(63,56){\makebox(0,0)[bl]{analysis of}}
\put(63,52){\makebox(0,0)[bl]{results}}


\put(00,10){\line(1,0){25}} \put(00,25){\line(1,0){25}}
\put(00,10){\line(0,1){15}} \put(25,10){\line(0,1){15}}

\put(12.5,25){\vector(0,1){23}}

\put(01,20){\makebox(0,0)[bl]{Problem}}
\put(01,16){\makebox(0,0)[bl]{reformulation/}}
\put(01,12){\makebox(0,0)[bl]{restructuring}}


\put(45,36.5){\oval(22,15)}

\put(36,39){\makebox(0,0)[bl]{Analysis of}}
\put(35,35){\makebox(0,0)[bl]{intermediate}}
\put(39.5,31){\makebox(0,0)[bl]{results}}

\put(45,29){\vector(0,-1){4}}

\put(35,32){\vector(-1,-1){10}}


\put(30,10){\line(1,0){26}} \put(30,25){\line(1,0){26}}
\put(30,10){\line(0,1){15}} \put(56,10){\line(0,1){15}}

\put(31,20){\makebox(0,0)[bl]{Reconfiguration}}
\put(31,16){\makebox(0,0)[bl]{of problem }}
\put(31,12){\makebox(0,0)[bl]{solving process}}

\put(32,25){\vector(0,1){23}}


\put(67,48){\line(0,-1){30.5}}

\put(67,17.5){\vector(-1,0){11}}


\put(75,48){\line(0,-1){42}}

\put(75,06){\line(-1,0){62.5}}

\put(12.5,06){\vector(0,1){4}}

\end{picture}
\end{center}

\subsection{History of DSS COMBI}

 A preliminary version of DSS COMBI was implemented as a set of
 multicriteria techniques (Fortran, mainframe, methods:
 several types of utility functions,
  Electre-like technique).
 Analysis, comparison  and aggregation of results, obtained via
 different techniques, was widely used.
 Further, DSS with method composition was designed.
 DSS COMBI was targeted to three type of resultant decisions (Fig. 23):
 (a) linear ranking
 (mainly, an intermediate result) \(\overline{B}\),
 (b) multicriteria ranking as sorting (ordinal priorities of alternatives,
 i.e., a layered structure) \(B\), and
 (c) multiticriteria ranking as a layered structure
  with intersection
  the layers (ordinal priorities of alternatives over an ordinal decision scale)
  \(\widetilde{B}\).

\begin{center}
\begin{picture}(149,88)

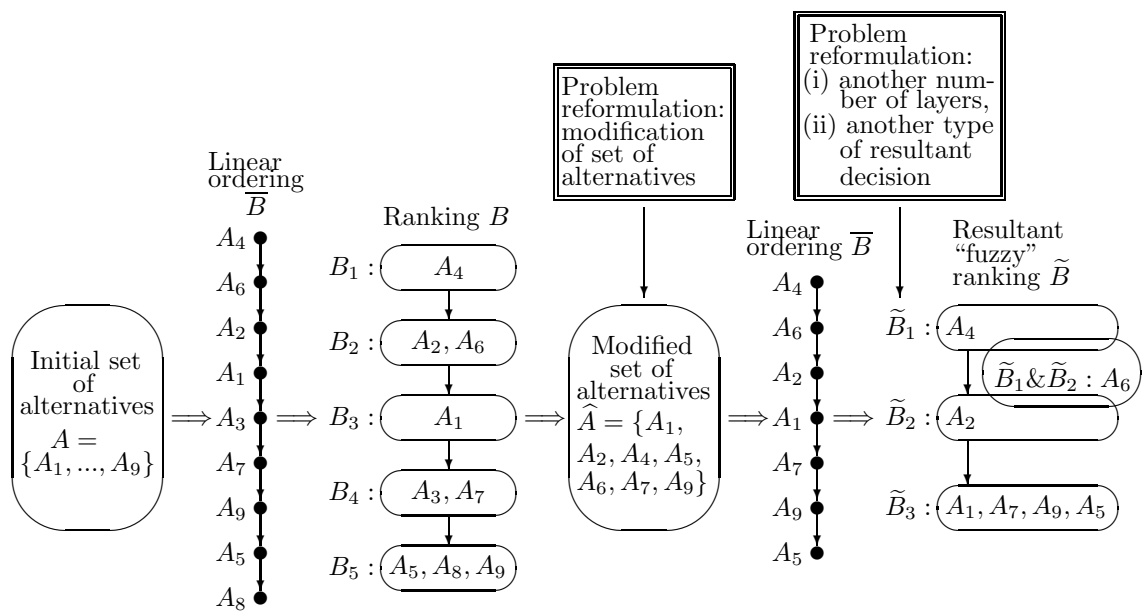
\put(30,00){\makebox(0,0)[bl]{Fig. 22. Series process with problem
 reformulation}}

\put(02.4,40.5){\makebox(0,0)[bl]{Initial set}}
\put(01,37.5){\makebox(0,0)[bl]{~~~~~ of}}

\put(01,34.5){\makebox(0,0)[bl]{alternatives}}

\put(05,30){\makebox(0,0)[bl]{\(A=\)}}

\put(01,26){\makebox(0,0)[bl]{\(\{A_{1},...,A_{9}\}\)}}

\put(10,34){\oval(20,30)}


\put(21,33){\makebox(0,0)[bl]{\(\Longrightarrow \)}}


\put(26,67){\makebox(0,0)[bl]{Linear}}

\put(26,64){\makebox(0,0)[bl]{ordering}}

\put(31,61){\makebox(0,0)[bl]{\(\overline{B}\)}}


\put(27,56.4){\makebox(0,0)[bl]{\(A_{4}\)}}

\put(33,58){\circle*{1.8}}


\put(33,58){\vector(0,-1){5}}

\put(27,50.4){\makebox(0,0)[bl]{\(A_{6}\)}}

\put(33,52){\circle*{1.8}}


\put(33,52){\vector(0,-1){5}}

\put(27,44.4){\makebox(0,0)[bl]{\(A_{2}\)}}

\put(33,46){\circle*{1.8}}


\put(33,46){\vector(0,-1){5}}

\put(27,38.4){\makebox(0,0)[bl]{\(A_{1}\)}}

\put(33,40){\circle*{1.8}}


\put(33,40){\vector(0,-1){5}}

\put(27,32.4){\makebox(0,0)[bl]{\(A_{3}\)}}

\put(33,34){\circle*{1.8}}


\put(33,34){\vector(0,-1){5}}

\put(27,26.4){\makebox(0,0)[bl]{\(A_{7}\)}}

\put(33,28){\circle*{1.8}}


\put(33,28){\vector(0,-1){5}}

\put(27,20.4){\makebox(0,0)[bl]{\(A_{9}\)}}

\put(33,22){\circle*{1.8}}


\put(33,22){\vector(0,-1){5}}

\put(27,14.4){\makebox(0,0)[bl]{\(A_{5}\)}}

\put(33,16){\circle*{1.8}}


\put(33,16){\vector(0,-1){5}}

\put(27,08.4){\makebox(0,0)[bl]{\(A_{8}\)}}

\put(33,10){\circle*{1.8}}



\put(35,33){\makebox(0,0)[bl]{\(\Longrightarrow \)}}


\put(49.3,59){\makebox(0,0)[bl]{Ranking \(B\)}}


\put(42,52.5){\makebox(0,0)[bl]{\(B_{1}:\)}}

\put(56,52.5){\makebox(0,0)[bl]{\(A_{4}\)}}

\put(58,54){\oval(18,6)}


\put(42,42.5){\makebox(0,0)[bl]{\(B_{2}:\)}}

\put(58,51){\vector(0,-1){4}}

\put(53,42.5){\makebox(0,0)[bl]{\(A_{2},A_{6}\)}}

\put(58,44){\oval(18,6)}


\put(42,32.5){\makebox(0,0)[bl]{\(B_{3}:\)}}

\put(58,41){\vector(0,-1){4}}

\put(56,32.5){\makebox(0,0)[bl]{\(A_{1}\)}}

\put(58,34){\oval(18,6)}


\put(42,22.5){\makebox(0,0)[bl]{\(B_{4}:\)}}

\put(58,31){\vector(0,-1){4}}

\put(53,22.5){\makebox(0,0)[bl]{\(A_{3},A_{7}\)}}

\put(58,24){\oval(18,6)}


\put(42,12.5){\makebox(0,0)[bl]{\(B_{5}:\)}}

\put(58,21){\vector(0,-1){4}}

\put(50,12.5){\makebox(0,0)[bl]{\(A_{5},A_{8},A_{9}\)}}

\put(58,14){\oval(18,6)}


\put(68,33){\makebox(0,0)[bl]{\(\Longrightarrow \)}}


\put(72,63){\line(1,0){24}} \put(72,81){\line(1,0){24}}
\put(72,63){\line(0,1){18}} \put(96,63){\line(0,1){18}}

\put(72.5,63.5){\line(1,0){23}} \put(72.5,80.5){\line(1,0){23}}
\put(72.5,63.5){\line(0,1){17}} \put(95.5,63.5){\line(0,1){17}}

\put(73,77){\makebox(0,0)[bl]{Problem}}
\put(73,74){\makebox(0,0)[bl]{reformulation:}}
\put(73,71){\makebox(0,0)[bl]{modification}}
\put(73,68){\makebox(0,0)[bl]{of set of}}
\put(73,65){\makebox(0,0)[bl]{alternatives}}

\put(84,62){\vector(0,-1){12}}


\put(104,63){\line(1,0){28}} \put(104,88){\line(1,0){28}}
\put(104,63){\line(0,1){25}} \put(132,63){\line(0,1){25}}

\put(104.5,63.5){\line(1,0){27}} \put(104.5,87.5){\line(1,0){27}}
\put(104.5,63.5){\line(0,1){24}} \put(131.5,63.5){\line(0,1){24}}

\put(106,84){\makebox(0,0)[bl]{Problem}}
\put(106,81){\makebox(0,0)[bl]{reformulation:}}
\put(105,77){\makebox(0,0)[bl]{(i) another num-}}
\put(110,74.5){\makebox(0,0)[bl]{ber of layers,}}
\put(105,71){\makebox(0,0)[bl]{(ii) another type }}
\put(110,68.5){\makebox(0,0)[bl]{of resultant}}
\put(110,65){\makebox(0,0)[bl]{decision}}

\put(118,62){\vector(0,-1){12}}


\put(77,42.5){\makebox(0,0)[bl]{Modified }}
\put(79.5,39.5){\makebox(0,0)[bl]{set of}}
\put(75,36.5){\makebox(0,0)[bl]{alternatives}}
\put(75,31.5){\makebox(0,0)[bl]{\(\widehat{A}=\{A_{1},\)}}
\put(75,27.5){\makebox(0,0)[bl]{\(A_{2},A_{4},A_{5},\)}}
\put(75,23.5){\makebox(0,0)[bl]{\(A_{6},A_{7},A_{9}\}\)}}

\put(84,34){\oval(20,30)}


\put(95,33){\makebox(0,0)[bl]{\(\Longrightarrow \)}}


\put(97.6,58){\makebox(0,0)[bl]{Linear}}
\put(97.6,55){\makebox(0,0)[bl]{ordering \(\overline{B}\)}}


\put(101,50.4){\makebox(0,0)[bl]{\(A_{4}\)}}

\put(107,52){\circle*{1.8}}


\put(107,52){\vector(0,-1){5}}

\put(101,44.4){\makebox(0,0)[bl]{\(A_{6}\)}}

\put(107,46){\circle*{1.8}}


\put(107,46){\vector(0,-1){5}}

\put(101,38.4){\makebox(0,0)[bl]{\(A_{2}\)}}

\put(107,40){\circle*{1.8}}


\put(107,40){\vector(0,-1){5}}

\put(101,32.4){\makebox(0,0)[bl]{\(A_{1}\)}}

\put(107,34){\circle*{1.8}}


\put(107,34){\vector(0,-1){5}}

\put(101,26.4){\makebox(0,0)[bl]{\(A_{7}\)}}

\put(107,28){\circle*{1.8}}


\put(107,28){\vector(0,-1){5}}

\put(101,20.4){\makebox(0,0)[bl]{\(A_{9}\)}}

\put(107,22){\circle*{1.8}}


\put(107,22){\vector(0,-1){5}}

\put(101,14.4){\makebox(0,0)[bl]{\(A_{5}\)}}

\put(107,16){\circle*{1.8}}


\put(109,33){\makebox(0,0)[bl]{\(\Longrightarrow \)}}



\put(125,58){\makebox(0,0)[bl]{Resultant}}
\put(125,54){\makebox(0,0)[bl]{``fuzzy''}}
\put(125,51){\makebox(0,0)[bl]{ranking \(\widetilde{B}\)}}


\put(116,44.5){\makebox(0,0)[bl]{\(\widetilde{B}_{1}:\)}}

\put(124,44.5){\makebox(0,0)[bl]{\(A_{4}\)}}

\put(135,46){\oval(24,6)}


\put(139.5,40){\oval(21,09)}

\put(130.5,37.5){\makebox(0,0)[bl]{\(\widetilde{B}_{1}\&\widetilde{B}_{2}:
A_{6}\)}}


\put(116,32.5){\makebox(0,0)[bl]{\(\widetilde{B}_{2}:\)}}

\put(127,43){\vector(0,-1){6}}

\put(124,32.5){\makebox(0,0)[bl]{\(A_{2}\)}}

\put(135,34){\oval(24,6)}


\put(116,20.5){\makebox(0,0)[bl]{\(\widetilde{B}_{3}:\)}}

\put(127,31){\vector(0,-1){6}}

\put(124,20.5){\makebox(0,0)[bl]{\(A_{1},A_{7},A_{9},A_{5}\)}}

\put(135,22){\oval(24,6)}

\end{picture}
\end{center}

\begin{center}
\begin{picture}(94,50)
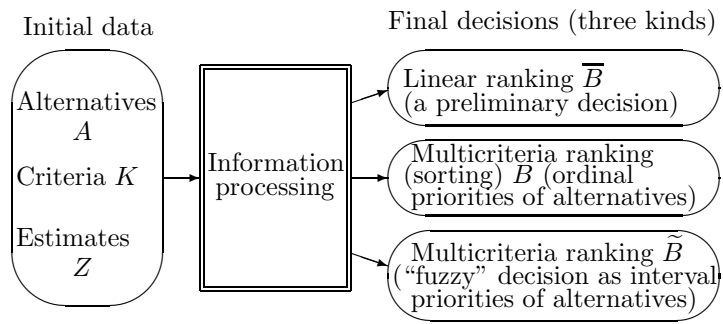
\put(00,00){\makebox(0,0)[bl] {Fig. 23. Solving framework of DSS
COMBI
 (\cite{lev86},\cite{lev98},\cite{levmih88})}}


\put(01.5,44){\makebox(0,0)[bl]{Initial data}}

\put(10,25){\oval(20,34)}

\put(0.5,34){\makebox(0,0)[bl]{Alternatives}}
\put(08,30){\makebox(0,0)[bl]{\(A\)}}

\put(0.5,24){\makebox(0,0)[bl]{Criteria \(K\)}}

\put(0.5,16){\makebox(0,0)[bl]{Estimates}}
\put(08,12){\makebox(0,0)[bl]{\(Z\)}}

\put(20,25){\vector(1,0){5}}


\put(25,10){\line(1,0){20}} \put(25,40){\line(1,0){20}}
\put(25,10){\line(0,1){30}} \put(45,10){\line(0,1){30}}

\put(25.5,10.5){\line(1,0){19}} \put(25.5,39.5){\line(1,0){19}}
\put(25.5,10.5){\line(0,1){29}} \put(44.5,10.5){\line(0,1){29}}

\put(26,26){\makebox(0,0)[bl]{Information}}
\put(27,22){\makebox(0,0)[bl]{processing}}

\put(45,15){\vector(3,-1){5}} \put(45,25){\vector(1,0){5}}
\put(45,35){\vector(3,1){5}}


\put(50,44){\makebox(0,0)[bl]{Final decisions (three kinds)}}


\put(72,37){\oval(44,10)}

\put(52,36.5){\makebox(0,0)[bl]{Linear ranking \(\overline{B}\)}}

\put(52,33){\makebox(0,0)[bl]{(a preliminary decision)}}


\put(72,25){\oval(44,10)}

\put(53,26.5){\makebox(0,0)[bl]{Multicriteria ranking}}
\put(52,23.5){\makebox(0,0)[bl]{(sorting) \(B\) (ordinal}}
\put(53,20.5){\makebox(0,0)[bl]{priorities of alternatives)}}


\put(72,12){\oval(44,12)}

\put(53,13.3){\makebox(0,0)[bl]{Multicriteria ranking
\(\widetilde{B}\)}}

\put(50.5,10){\makebox(0,0)[bl]{(``fuzzy'' decision as interval}}

\put(53,07){\makebox(0,0)[bl]{priorities of alternatives) }}

\end{picture}
\end{center}

 Table 3 presents a brief description of DSS COMBI generations
 with method composition
 (\cite{lev93},\cite{lev98},\cite{levmih88}).

 A functional graph menu was realized in
 DSS COMBI (since generation 1) (Fig. 24).

\begin{center}
\begin{picture}(144.5,52)
\put(40,48){\makebox(0,0)[bl]{Table 3. Generations of DSS COMBI}}

\put(00,00){\line(1,0){144.5}} \put(00,34){\line(1,0){144.5}}
\put(00,46){\line(1,0){144.5}}

\put(00,0){\line(0,1){46}} \put(25,0){\line(0,1){46}}
\put(42,0){\line(0,1){46}} \put(60,0){\line(0,1){46}}
\put(74,0){\line(0,1){46}} \put(84,0){\line(0,1){46}}
\put(93,0){\line(0,1){46}} \put(110,0){\line(0,1){46}}
\put(128.5,0){\line(0,1){46}} \put(144.5,0){\line(0,1){46}}


\put(01,42){\makebox(0,0)[bl]{Generation of}}
\put(01,39){\makebox(0,0)[bl]{DSS COMBI}}

\put(26,42){\makebox(0,0)[bl]{Type of}}
\put(26,39){\makebox(0,0)[bl]{computer}}

\put(43,42){\makebox(0,0)[bl]{Type of }}
\put(43,39.5){\makebox(0,0)[bl]{interface}}

\put(61,42){\makebox(0,0)[bl]{Domain}}

\put(75,42){\makebox(0,0)[bl]{Lear-}}
\put(75,39){\makebox(0,0)[bl]{ning }}

\put(85,42){\makebox(0,0)[bl]{Year}}

\put(94,42){\makebox(0,0)[bl]{Reference}}

\put(111,42){\makebox(0,0)[bl]{Presenta-}}
\put(111,39){\makebox(0,0)[bl]{tion at}}
\put(111,36){\makebox(0,0)[bl]{conference}}

\put(129,42){\makebox(0,0)[bl]{Usage in }}
\put(129,39){\makebox(0,0)[bl]{teaching}}


\put(01,30){\makebox(0,0)[bl]{0. COMBI}}
\put(02,26){\makebox(0,0)[bl]{(Pascal-based)}}

\put(26,30){\makebox(0,0)[bl]{Mini-}}
\put(26,26){\makebox(0,0)[bl]{computer}}

\put(43,30){\makebox(0,0)[bl]{Language-}}
\put(43,26.7){\makebox(0,0)[bl]{based}}

\put(61,30){\makebox(0,0)[bl]{Various}}

\put(75,30){\makebox(0,0)[bl]{None}}

\put(85,30){\makebox(0,0)[bl]{1987}}

\put(93.5,29.3){\makebox(0,0)[bl]{\cite{levmih87a},\cite{levmih87b}}}
\put(93.5,25.3){\makebox(0,0)[bl]{\cite{levmih90}}}

\put(110.5,30){\makebox(0,0)[bl]{None}}

\put(129,30){\makebox(0,0)[bl]{None}}


\put(01,22){\makebox(0,0)[bl]{1. COMBI PC }}
\put(02,18){\makebox(0,0)[bl]{(Pascal-based)}}

\put(26,22){\makebox(0,0)[bl]{PC}}

\put(43,22){\makebox(0,0)[bl]{Graphical}}
\put(43,18.7){\makebox(0,0)[bl]{menu}}

\put(61,22){\makebox(0,0)[bl]{Various}}

\put(75,22){\makebox(0,0)[bl]{Yes}}

\put(85,22){\makebox(0,0)[bl]{1988}}

\put(93.5,21.3){\makebox(0,0)[bl]{\cite{lev86},\cite{levmih88}}}

\put(110.5,22){\makebox(0,0)[bl]{SPUDM-89}}

\put(129,22){\makebox(0,0)[bl]{None}}


\put(01,14){\makebox(0,0)[bl]{2. COMBI PC}}
\put(02,10){\makebox(0,0)[bl]{(C-based)}}

\put(26,14){\makebox(0,0)[bl]{PC}}

\put(43,14){\makebox(0,0)[bl]{Graphical}}
\put(43,10.7){\makebox(0,0)[bl]{menu}}

\put(61,14){\makebox(0,0)[bl]{Various}}

\put(75,14){\makebox(0,0)[bl]{Yes}}

\put(85,14){\makebox(0,0)[bl]{1989}}

\put(93.5,13.3){\makebox(0,0)[bl]{\cite{lev93},\cite{lev94},}}
\put(93.5,09.3){\makebox(0,0)[bl]{\cite{lev98},\cite{levmih88}}}

\put(110.5,14){\makebox(0,0)[bl]{MCDM-90}}
\put(110.5,10){\makebox(0,0)[bl]{EWHCI-93}}

\put(129,13.3){\makebox(0,0)[bl]{\cite{lev95edu},\cite{lev96edu}}}
\put(129,09.3){\makebox(0,0)[bl]{\cite{lev06edu},\cite{lev10edu}}}


\put(01,06){\makebox(0,0)[bl]{3. COMBI PC}}
\put(02,02){\makebox(0,0)[bl]{(C-based)}}

\put(26,06){\makebox(0,0)[bl]{PC}}

\put(43,06){\makebox(0,0)[bl]{Graphical}}
\put(43,02.7){\makebox(0,0)[bl]{menu}}

\put(61,06){\makebox(0,0)[bl]{Invest-}}
\put(61,02){\makebox(0,0)[bl]{ment}}

\put(75,06){\makebox(0,0)[bl]{Yes}}

\put(85,06){\makebox(0,0)[bl]{1991}}

\put(93.5,05.3){\makebox(0,0)[bl]{\cite{lev93}}}

\put(110.5,06){\makebox(0,0)[bl]{None}}

\put(129,06){\makebox(0,0)[bl]{None}}

\end{picture}
\end{center}

\begin{center}
\begin{picture}(108,80)
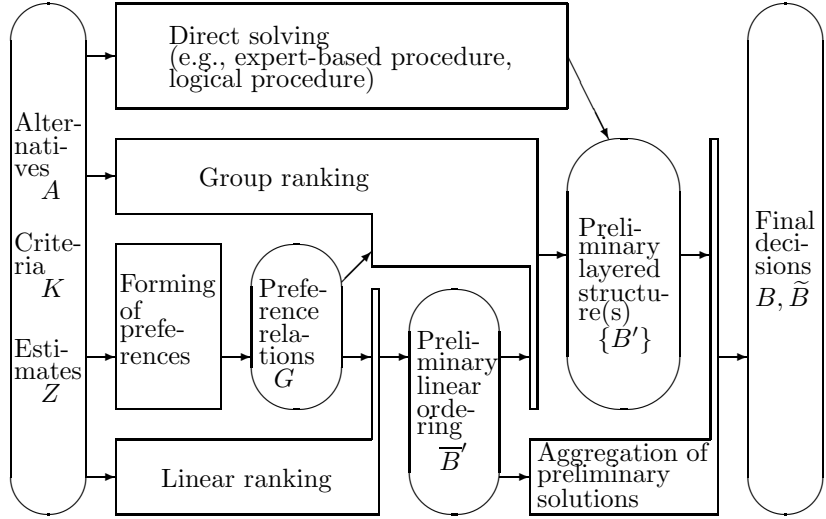
\put(01,00){\makebox(0,0)[bl] {Fig. 24. Structure of functional
 planar graph menu (\cite{lev86},\cite{lev98},\cite{levmih88})}}


\put(05,41){\oval(10,68)}

\put(0.5,58){\makebox(0,0)[bl]{Alter-}}
\put(0.5,55){\makebox(0,0)[bl]{nati-}}
\put(0.5,52){\makebox(0,0)[bl]{ves}}
\put(04,49){\makebox(0,0)[bl]{\(A\)}}

\put(0.5,42){\makebox(0,0)[bl]{Crite-}}
\put(0.5,39){\makebox(0,0)[bl]{ria}}
\put(04,36){\makebox(0,0)[bl]{\(K\)}}

\put(0.5,28){\makebox(0,0)[bl]{Esti-}}
\put(0.5,25){\makebox(0,0)[bl]{mates}}
\put(04,22){\makebox(0,0)[bl]{\(Z\)}}


\put(103,41){\oval(10,68)}

\put(99,45){\makebox(0,0)[bl]{Final}}
\put(99,42){\makebox(0,0)[bl]{deci-}}
\put(99,39){\makebox(0,0)[bl]{sions}}
\put(99,34){\makebox(0,0)[bl]{\(B,\widetilde{B}\)}}

\put(94,28){\vector(1,0){4}}


\put(48,37){\line(1,0){1}}

\put(48,17){\line(0,1){20}}

\put(14,07){\line(0,1){10}} \put(49,07){\line(0,1){30}}

\put(14,17){\line(1,0){34}} \put(14,07){\line(1,0){35}}

\put(20,10){\makebox(0,0)[bl]{Linear ranking}}

\put(10,12){\vector(1,0){4}}


\put(14,21){\line(0,1){22}} \put(28,21){\line(0,1){22}}
\put(14,21){\line(1,0){14}} \put(14,43){\line(1,0){14}}

\put(14.5,36){\makebox(0,0)[bl]{Forming}}
\put(14.5,33){\makebox(0,0)[bl]{of}}
\put(14.5,30){\makebox(0,0)[bl]{prefe-}}
\put(14.5,27){\makebox(0,0)[bl]{rences}}

\put(10,28){\vector(1,0){4}}

\put(28,28){\vector(1,0){4}}

\put(44,28){\vector(1,0){4}}

\put(49,28){\vector(1,0){4}}


\put(38,32){\oval(12,22)}

\put(33,36){\makebox(0,0)[bl]{Prefe-}}
\put(33,33){\makebox(0,0)[bl]{rence}}
\put(33,30){\makebox(0,0)[bl]{rela-}}
\put(33,27){\makebox(0,0)[bl]{tions}}
\put(35,24){\makebox(0,0)[bl]{\(G\)}}


\put(59,22){\oval(12,30)}

\put(54,29){\makebox(0,0)[bl]{Preli-}}
\put(54,25.6){\makebox(0,0)[bl]{minary}}
\put(54,23){\makebox(0,0)[bl]{linear}}
\put(54,20){\makebox(0,0)[bl]{orde-}}
\put(54,17){\makebox(0,0)[bl]{ring}}
\put(57,13){\makebox(0,0)[bl]{\(\overline{B}'\)}}

\put(65,28){\vector(1,0){4}}


\put(93,57){\line(1,0){1}}

\put(93,17){\line(0,1){40}}

\put(69,07){\line(0,1){10}} \put(94,07){\line(0,1){50}}

\put(69,17){\line(1,0){24}} \put(69,07){\line(1,0){25}}

\put(70,13.5){\makebox(0,0)[bl]{Aggregation of}}
\put(70,10.6){\makebox(0,0)[bl]{preliminary}}
\put(70,08){\makebox(0,0)[bl]{solutions}}

\put(65,12){\vector(1,0){4}}


\put(14,47){\line(0,1){10}}

\put(69,40){\line(-1,0){21}}

\put(48,40){\line(0,1){7}}

\put(69,21){\line(0,1){19}} \put(70,21){\line(0,1){36}}
\put(69,21){\line(1,0){1}}

\put(14,47){\line(1,0){34}} \put(14,57){\line(1,0){56}}

\put(25,50){\makebox(0,0)[bl]{Group ranking}}

\put(10,52){\vector(1,0){4}}

\put(44,38){\vector(1,1){4}}

\put(70,41.5){\vector(1,0){4}}


\put(14,61){\line(0,1){14}} \put(74,61){\line(0,1){14}}

\put(14,61){\line(1,0){60}} \put(14,75){\line(1,0){60}}

\put(21,69){\makebox(0,0)[bl]{Direct solving}}
\put(21,66){\makebox(0,0)[bl]{(e.g., expert-based procedure,}}
\put(21,63){\makebox(0,0)[bl]{logical procedure)}}

\put(10,68){\vector(1,0){4}} \put(74,68){\vector(1,-2){5.5}}


\put(81.5,39){\oval(15,36)}

\put(75.5,44){\makebox(0,0)[bl]{Preli-}}
\put(75.5,41){\makebox(0,0)[bl]{minary}}
\put(75.5,38){\makebox(0,0)[bl]{layered}}
\put(75.5,35.4){\makebox(0,0)[bl]{structu-}}
\put(75.5,32){\makebox(0,0)[bl]{re(s)}}
\put(78,28.5){\makebox(0,0)[bl]{\(\{B'\}\)}}

\put(89,41.5){\vector(1,0){4}}

\end{picture}
\end{center}

\subsection{Main Components of DSS}

 Here, the following DSS parts are considered:
 {\it 1.} information part as data, knowledge;
 {\it 2.} operational part as tools for information processing);
 and
 {\it 3.} human part as user or group of users
 (including experts, decision maker, etc.).

 Usually information part includes the following:
   (1) data (alternatives or basic items, criteria, multicriteria
       estimates of alternatives upon criteria, preference relations);
   (2) tools for maintaining data
       (DBMS and interfaces with other commercial DBMSs);
   (3) support information for learning
       (e.g., a helper, etc.).

 In this paper (as in DSS COMBI),
 the following series (or series-parallel) framework
 of information (a preference relation or a matrix)
  processing is examined
 (\cite{lev93}, \cite{lev98}, \cite{levmih88}):
 ~(1) basic data as alternatives (items) (\(A\)), criteria (\(K\)), multicriteria
        estimates of alternatives upon criteria (\(Z\));
 ~(2) preference relation of alternatives (\(G\));
 ~(3) an intermediate linear ordering of alternatives (\(\overline{B}\))
 ~(4) an preliminary ranking of alternatives as a layered structure(s) (\(B'\));
 and
 ~(5) resultant ranking of alternatives
    as layered structure
    (\(B\))
     (an ordinal priority for each alternative)
        or
        fuzzy ranking  (\(\widetilde{B}\))
        (with intersection of layers).

 The following kinds of basic operations are considered:
   (i) data processing (series and/or parallel);
   (ii) data aggregation.
 In addition,
    parallelization of the
         solving process on the basis of various components
         (alternatives, criteria, techniques, experts) can be used.
 The solving process may be presented as a hierarchy with the following
 functional/operatiopnal layers (\cite{lev93}, \cite{lev98}, \cite{levmih88}):
   (1) algorithms and man-machine interactive procedures for data
       transformation (bottom layer);
   (2) strategies (step-by-step schemes of data transformation,
       particularly series-parallel ones);
   (3) scenarios (complexes of strategies with their analysis and feedback).

 Thus, the following alternative series frameworks
 (composite solving strategies)
 of information processing are considered (Fig. 11, Fig. 24):
%

   {\bf 1.} Basic series framework:~
   \(E: ~  A \Rightarrow G \Rightarrow \overline{B}' \Rightarrow \{B'\} \Rightarrow B\) or \( \widetilde{B}
   \).

  {\bf 2.} Compressed series frameworks:

 ~~~{\it 2.1.}~
 \( W^{1}:~  A \Rightarrow \{ B'\} \Rightarrow B\) or \(\widetilde{B}\).

 ~~~{\it 2.2.}~
 \( W^{2}:~  A \Rightarrow G \Rightarrow \{ B'\} \Rightarrow B\) or \( \widetilde{B}\).

 ~~~{\it 2.3.}~
   \( W^{3}:~ A \Rightarrow G \Rightarrow \overline{B}' \Rightarrow B\) or \( \widetilde{B}\).

 ~~~{\it 2.4.}~
   \( W^{4}:~ A \Rightarrow  \overline{B}' \Rightarrow \{ B'\} \Rightarrow B\) or \(  \widetilde{B}
   \).

  ~~~{\it 2.5.}~
   \( W^{5}:~  A \Rightarrow \overline{B}' \Rightarrow \{ B'\} \Rightarrow B\) or \(   \widetilde{B}
   \).

  {\bf 3.} Direct solving process:~ \( D:~ A \Rightarrow B\) or \(\widetilde{B}\).

 Thus, the solving process can be considered as the following
 four-part system (Fig. 25):
  ~~\( S = H \star T \star U \star X\),~
  where the following four stages are basic ones:
 ~(i)  stage \(H\) corresponds to forming a preliminary preference relation
 \(G\)
  (over alternatives \(A\))
 (an algorithm or a procedure);
  ~(ii) stage \(T\) corresponds to forming a preliminary linear ranking \(\overline{B}\)
 (an algorithm or a procedure);
 ~(iii) stage \(U\) corresponds to forming some preliminary  rankings \( \{B \}\)
 (an algorithm or a procedure); and
 ~(iv) stage \(X\) corresponds to aggregation  of the preliminary  rankings \( \{B \}\)
 into the resultant decisions
 \(B\) (or \(\widetilde{B}\))
  (an algorithm or a procedure).
%
 The following basic local techniques (as processing units) have been
 used
  in DSS COMBI (\cite{lev93}, \cite{lev98}, \cite{levmih88}):

 {\bf I.} Stage \(H\):~
  absent (\(H_{0}\)),
  pairwise comparison (a simple expert-based procedure) (\(H_{1}\)),
  dominance by  of Pareto-rule  (\(H_{2}\))
   (e.g., \cite{mirkin79},\cite{pareto71}),
  outranking Electre-like technique
 (a special Electre-based interactive procedure with feedback, designed by M.Sh. Levin
 \cite{levmih88})
  (\(H_{3}\)).

{\bf II.} Stage \(T\):~
 absent (\(T_{0}\)),
 line elements sum of preference matrix (\(T_{1}\)),
  additive utility function
 (e.g., \cite{fis70},\cite{kee76})
 (\(T_{2}\)).

 {\bf III.} Stage \(U\):~
  absent (\(U_{0}\)),
  step-by-step revelation of {\it maximal} elements (\(U_{1}\)),
 step-by-step revelation of Pareto-efficient elements
 (e.g., \cite{mirkin79},\cite{pareto71})
  (\(U_{2}\)),
 dividing the linear ranking (\(U_{3}\)),
 expert procedure for ranking
 (expert-based direct solving procedure
 (e.g., \cite{larmf86})
 (\(U_{4}\)),
 direct logical method for ranking
 (direct solving method based on logical approach, designed by A.A.
 Mikhailov)
 (\cite{levmih87c},\cite{levmih88})
 (\(U_{5}\)).

{\bf IV.} Stage \(X\):~
 absent (\(X_{0}\)),
 simple election-like procedure (\(X_{1}\)),
 aggregation model based on special knapsack-like problem (\(X_{2}\))
 \cite{lev98}.

\begin{center}
\begin{picture}(70,42)
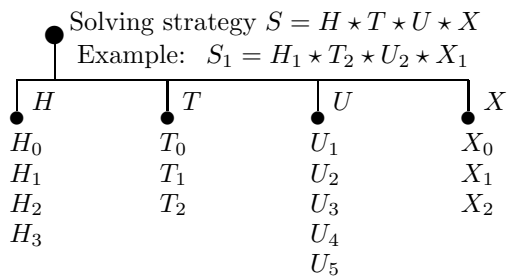
\put(02,00){\makebox(0,0)[bl]{Fig. 25. Morphology
 of solving strategy}}


\put(1,28){\line(0,1){4}} \put(01,27){\circle*{1.6}}
\put(03,28){\makebox(0,0)[bl]{\(H\)}}

\put(00,22){\makebox(0,0)[bl]{\(H_{0}\)}}
\put(00,18){\makebox(0,0)[bl]{\(H_{1}\)}}
\put(00,14){\makebox(0,0)[bl]{\(H_{2}\)}}
\put(00,10){\makebox(0,0)[bl]{\(H_{3}\)}}


\put(21,28){\line(0,1){4}} \put(21,27){\circle*{1.6}}
\put(23,28){\makebox(0,0)[bl]{\(T\)}}

\put(20,22){\makebox(0,0)[bl]{\(T_{0}\)}}
\put(20,18){\makebox(0,0)[bl]{\(T_{1}\)}}
\put(20,14){\makebox(0,0)[bl]{\(T_{2}\)}}


\put(41,28){\line(0,1){4}} \put(41,27){\circle*{1.6}}
\put(43,28){\makebox(0,0)[bl]{\(U\)}}

\put(40,22){\makebox(0,0)[bl]{\(U_{1}\)}}
\put(40,18){\makebox(0,0)[bl]{\(U_{2}\)}}
\put(40,14){\makebox(0,0)[bl]{\(U_{3}\)}}
\put(40,10){\makebox(0,0)[bl]{\(U_{4}\)}}
\put(40,06){\makebox(0,0)[bl]{\(U_{5}\)}}


\put(61,28){\line(0,1){4}} \put(61,27){\circle*{1.6}}
\put(63,28){\makebox(0,0)[bl]{\(X\)}}

\put(60,22){\makebox(0,0)[bl]{\(X_{0}\)}}
\put(60,18){\makebox(0,0)[bl]{\(X_{1}\)}}
\put(60,14){\makebox(0,0)[bl]{\(X_{2}\)}}


\put(01,32){\line(1,0){60}} \put(06,32){\line(0,1){06}}
\put(06,38){\circle*{2.5}}

\put(08,38){\makebox(0,0)[bl]{Solving strategy
 \(S = H \star T \star U \star X \)}}

\put(09,33.6){\makebox(0,0)[bl] {Example: ~\(S_{1}= H_{1} \star
T_{2} \star U_{2} \star X_{1} \)}}

\end{picture}
\end{center}

 Here the index \(0\) corresponds to absence of the precessing at
 the stage.
 Thus, the following examples can be considered
 (symbol \(\&\) corresponds to parallel integration of processing operations):

 (1) a series solving  strategy
 ~\(S_{1}= H_{1} \star T_{2} \star U_{2} \star X_{0} \);

 (2) a series-parallel solving strategy:

 \(S_{2} =
 ( S' \& S'' \& S''') \star X_{1} =\)
 \( ( ( H_{3} \star T_{1} \star U_{1}) \&
 ( H_{3} \star T_{1} \star U_{2}) \&
 (H_{1} \star T_{2} \star U_{1}) ) \star X_{1}\).

 Parallelization of the solving process is based on three
 approaches:

 (1) concurrent usage of different experts in
 the same interactive procedure (e.g., in \(U_{4}\));

 (2) concurrent usage of different methods at the same
 solving stage
  (e.g., \(H_{3}\) );

 (3) concurrent usage of the same method with different parameters
 at the same solving stage
  (e.g., in  \(U_{4}\)).

 The design of solving strategies consists of two problems \cite{lev98}:
 (a) selection of local techniques and
 (b) combinatorial synthesis of the selected techniques.
 The requirements to the composite solving strategies are based on the following:
  (i) kinds of task and of users;
  (ii) available resources (e.g., human, computer, time);
  (iii) features of the decision situation (e.g., kind of uncertainty,
      required precision and robustness of result).
 Thus, the following six criteria are examined \cite{lev98}:
 (1) required computer resources;
 (2) required human resources;
 (3) quality of ranking (robustness, etc.);
 (4) possibility for data representation;
 (5) possibility for an analysis of intermediate data; and
 (6) usability (easy to learn, understanding, etc.).

 Two direct solving strategies are as follows:
 (a) direct expert based ranking (e.g., \cite{larmf86}):~
  \(D_{1}  = H_{4} \star T_{5} \star U_{4} \star X_{0}\);
 (b)  direct logical method for ranking
 (suggested and designed by A.A. Mikhailov
 \cite{levmih88}):~
  \(D'_{2}  = H_{4} \star T_{6} \star U_{5} \star X_{0}\)
  or
 for the case of several experts:
 (a) \(D''_{2}  = H_{4} \star T_{6} \star U_{5} \star X_{1}\),
 (b) \(D'''_{2}  = H_{4} \star T_{6} \star U_{5} \star X_{1}\).
 In our design approach, the strategies are based on defined
 top-level compatibility between their components
 (and zero compatibility between the components of the strategies and
 other components).

 Note, some traditional approaches
 to build a solving process for decision making
  are oriented to the selection of the best method
 (or algorithm, model)
 (e.g., \cite{ban93},\cite{gui98},\cite{oz92},\cite{sho88}).
 A model composition on the basis of a filter space is described in
 \cite{cha02}.
 A non-linear recursive process consisting of four steps
 is analyzed for multicriteria decision aid in \cite{gui98}.

 In general, human part consists of the following:
   (a) user or group of users (experts, decision makers);
   (b) techniques for the modeling, diagnostics, selection, and
       assignment of specialists;
   (c) subsystem for user training;
   (d) user interfaces and tools for their adaptation.
 In recent years, many authors have been investigated user modeling
 (e.g., \cite{bar90},\cite{kob89},\cite{lev93},\cite{rou92}).
 Issues of adaptation of user interfaces are considered in
 (e.g., \cite{deb94},\cite{dos89},\cite{lev02},\cite{nor89},\cite{sch93}).
 Graphical interaction in multicriteria decision making is considered in
 (\cite{jon93},\cite{kor91},\cite{mar88},\cite{pra90},\cite{vet94}).
 The effectiveness of different representations for managerial problem
 solving has been studied in \cite{sme97}. Note a special prospective
 class of graph-based modeling systems for decision making processes
 has been proposed in (\cite{jon91},\cite{jon93},\cite{lee93}).
 A hierarchical process of the user interface design for DSS COMBI
 on the basis of HMMD, and comparing some system
 versions are described in details in (\cite{lev94},\cite{lev98}).

\section{SYNTHESIS OF COMPOSITE STRATEGY}

\subsection{Combinatorial Synthesis with Interval Multiset
 estimates}

 Interval multiset estimates have been suggested by M.Sh. Levin
 in \cite{lev12a}.
 The approach consists in assignment of elements (\(1,2,3,...\))
 into an ordinal scale \([1,2,...,l]\).
 As a result, a multi-set based estimate is obtained,
 where a basis set involves all levels of the ordinal scale:
 \(\Omega = \{ 1,2,...,l\}\) (the levels are linear ordered:
 \(1 \succ 2 \succ 3 \succ ...\)) and
 the assessment problem (for each alternative)
 consists in selection of a multiset over set \(\Omega\) while taking into
 account two conditions:

 {\it 1.} cardinality of the selected multiset equals a specified
 number of elements \( \eta = 1,2,3,...\)
 (i.e., multisets of cardinality \(\eta \) are considered);

 {\it 2.} ``configuration'' of the multiset is the following:
 the selected elements of \(\Omega\) cover an interval over scale \([1,l]\)
 (i.e., ``interval multiset estimate'').

 Thus, an estimate \(e\) for an alternative \(A\) is
 (scale \([1,l]\), position-based form or position form):
 \(e(A) = (\eta_{1},...,\eta_{\iota},...,\eta_{l})\),
 where \(\eta_{\iota}\) corresponds to the number of elements at the
 level \(\iota\) (\(\iota = \overline{1,l}\)), or
 \(e(A) = \{ \overbrace{1,...,1}^{\eta_{1}},\overbrace{2,...2}^{\eta_{2}},
 \overbrace{3,...,3}^{\eta_{3}},...,\overbrace{l,...,l}^{\eta_{l}}
 \}\).
 The number of multisets of cardinality \(\eta\),
 with elements taken from a finite set of cardinality \(l\),
 is called the
 ``multiset coefficient'' or ``multiset number''
  (\cite{knuth98},\cite{yager86}):
 ~~\( \mu^{l,\eta} =
   \frac{l(l+1)(l+2)... (l+\eta-1) } {\eta!}
   \).
 This number corresponds to possible estimates
 (without taking into account interval condition 2).
 In the case of condition 2,
 the number of estimates is decreased.
 Generally, assessment problems based on interval multiset estimates
 can be denoted as follows: ~\(P^{l,\eta}\).
 A poset-like scale of interval multiset estimates for assessment problem \(P^{3,4}\)
 is presented
 in Fig. 26.
 The assessment problem will be used in our applied numerical
 examples.

 In addition, operations over multiset estimates
 are used \cite{lev12a}:
 integration, vector-like proximity, aggregation, and alignment.


 Integration of estimates (mainly, for composite systems)
 is based on summarization of the estimates by components (i.e.,
 positions).
 Let us consider \(n\) estimates (position form):~~
 estimate \(e^{1} = (\eta^{1}_{1},...,\eta^{1}_{\iota},...,\eta^{1}_{l})
 \),
  {\bf . . .},
 estimate \(e^{\kappa} = (\eta^{\kappa}_{1},...,\eta^{\kappa}_{\iota},...,\eta^{\kappa}_{l})
 \),
  {\bf . . .},
 estimate \(e^{n} = (\eta^{n}_{1},...,\eta^{n}_{\iota},...,\eta^{n}_{l})
 \).
 Then, the integrated estimate is:~
 estimate \(e^{I} = (\eta^{I}_{1},...,\eta^{I}_{\iota},...,\eta^{I}_{l})
 \),
 where
 \(\eta^{I}_{\iota} = \sum_{\kappa=1}^{n} \eta^{\kappa}_{\iota} ~~ \forall
 \iota = \overline{1,l}\).
 In fact, the operation \(\biguplus\) is used for multiset estimates:
 \(e^{I} = e^{1} \biguplus ... \biguplus e^{\kappa} \biguplus ... \biguplus e^{n}\).


 Further, vector-like proximity is described.
  Let \(A_{1}\) and \(A_{2}\) be two alternatives
 with corresponding
 interval multiset estimates
 \(e(A_{1})\), \(e(A_{2})\).
  Vector-like proximity for the alternatives above is:
 ~~\(\delta ( e(A_{1}), e(A_{2})) = (\delta^{-}(A_{1},A_{2}),\delta^{+}(A_{1},A_{2}))\),
 where vector components are:
 (i) \(\delta^{-}\) is the number of one-step changes:
 element of quality \(\iota + 1\) into element of quality \(\iota\) (\(\iota = \overline{1,l-1}\))
 (this corresponds to ``improvement'');
 (ii) \(\delta^{+}\) is the number of one-step changes:
 element of quality \(\iota\) into element of quality  \(\iota+1\) (\(\iota = \overline{1,l-1}\))
 (this corresponds to ``degradation'').
 It is assumed:
 ~\( | \delta ( e(A_{1}), e(A_{2})) | = | \delta^{-}(A_{1},A_{2}) | + |\delta^{+}(A_{1},A_{2})|
 \).
%


 Now let us consider median estimates (aggregation)
  for the
 specified set of initial estimates
 (traditional approach).
 Let \(E = \{ e_{1},...,e_{\kappa},...,e_{n}\}\)
 be the set of specified estimates
 (or a corresponding set of specified alternatives),
 let \(D \)
 be the set of all possible estimates
 (or a corresponding set of possible alternatives)
 (\( E  \subseteq D \)).
  Thus, the median estimates
  (``generalized median'' \(M^{g}\) and ``set median'' \(M^{s}\)) are:
 ~~\(M^{g} =   \arg \min_{M \in D}~
   \sum_{\kappa=1}^{n} ~  | \delta (M, e_{\kappa}) |; ~~
 M^{s} =   \arg \min_{M\in E} ~
   \sum_{\kappa=1}^{n} ~ | \delta (M, e_{\kappa}) |\).

 A brief description  of combinatorial synthesis
 (HMMD) was presented in introduction.

 Let \(S\) be a system consisting of \(m\) parts (components):
 \(R(1),...,R(i),...,R(m)\).
 A set of design alternatives
 is generated for each system part above.

 The problem is:

~~

 {\it Find a composite design alternative}
 ~~ \(S=S(1)\star ...\star S(i)\star ...\star S(m)\)~~
 {\it of DAs (one representative design alternative}
 ~\(S(i)\)
 {\it for each system component/part}
  ~\(R(i)\), \(i=\overline{1,m}\)
  {\it )}
 {\it with non-zero compatibility}
 {\it between design alternatives.}

~~

 A discrete ``space'' of the system excellence
 (a poset)
 on the basis of the following vector is used:
 ~~\(N(S)=(w(S);e(S))\),
 ~where \(w(S)\) is the minimum of pairwise compatibility
 between DAs which correspond to different system components
 (i.e.,
 \(~\forall ~R_{j_{1}}\) and \( R_{j_{2}}\),
 \(1 \leq j_{1} \neq j_{2} \leq m\))
 in \(S\),
 ~\(e(S)=(\eta_{1},...,\eta_{\iota},...,\eta_{l})\),
 ~where \(\eta_{\iota}\) is the number of DAs of the \(\iota\)th quality in \(S\).

\begin{center}
\begin{picture}(70,139)

\put(04.6,00){\makebox(0,0)[bl] {Fig. 26. Scale, estimates
 (\(P^{3,4}\)) \cite{lev12a}}}


\put(25,126.7){\makebox(0,0)[bl]{\(e^{3,4}_{1}\) }}

\put(28,129){\oval(16,5)} \put(28,129){\oval(16.5,5.5)}


\put(49,127){\makebox(0,0)[bl]{\((4,0,0)\) }}

\put(00,128.5){\line(0,1){08}} \put(04,128.5){\line(0,1){08}}

\put(00,130.5){\line(1,0){4}} \put(00,132.5){\line(1,0){4}}
\put(00,134.5){\line(1,0){4}} \put(00,136.5){\line(1,0){4}}

\put(00,128.5){\line(1,0){12}}

\put(00,127){\line(0,1){3}} \put(04,127){\line(0,1){3}}
\put(08,127){\line(0,1){3}} \put(12,127){\line(0,1){3}}

\put(01.5,124.5){\makebox(0,0)[bl]{\(1\)}}
\put(05.5,124.5){\makebox(0,0)[bl]{\(2\)}}
\put(09.5,124.5){\makebox(0,0)[bl]{\(3\)}}


\put(28,120){\line(0,1){6}}

\put(25,114.7){\makebox(0,0)[bl]{\(e^{3,4}_{2}\) }}

\put(28,117){\oval(16,5)}


\put(49,115){\makebox(0,0)[bl]{\((3,1,0)\) }}

\put(00,116.5){\line(0,1){06}} \put(04,116.5){\line(0,1){06}}
\put(08,116.5){\line(0,1){02}}

\put(00,118.5){\line(1,0){8}} \put(00,120.5){\line(1,0){4}}
\put(00,122.5){\line(1,0){4}}

\put(00,116.5){\line(1,0){12}}

\put(00,115){\line(0,1){3}} \put(04,115){\line(0,1){3}}
\put(08,115){\line(0,1){3}} \put(12,115){\line(0,1){3}}

\put(01.5,112.5){\makebox(0,0)[bl]{\(1\)}}
\put(05.5,112.5){\makebox(0,0)[bl]{\(2\)}}
\put(09.5,112.5){\makebox(0,0)[bl]{\(3\)}}


\put(28,108){\line(0,1){6}}

\put(25,102.7){\makebox(0,0)[bl]{\(e^{3,4}_{3}\) }}

\put(28,105){\oval(16,5)}


\put(49,102){\makebox(0,0)[bl]{\((2,2,0)\) }}

\put(00,104.5){\line(0,1){04}}

\put(04,104.5){\line(0,1){04}} \put(08,104.5){\line(0,1){04}}
\put(00,106.5){\line(1,0){8}} \put(00,108.5){\line(1,0){8}}

\put(00,104.5){\line(1,0){12}}

\put(00,103){\line(0,1){3}} \put(04,103){\line(0,1){3}}
\put(08,103){\line(0,1){3}} \put(12,103){\line(0,1){3}}

\put(01.5,100.5){\makebox(0,0)[bl]{\(1\)}}
\put(05.5,100.5){\makebox(0,0)[bl]{\(2\)}}
\put(09.5,100.5){\makebox(0,0)[bl]{\(3\)}}


\put(28,96){\line(0,1){6}}

\put(25,90.7){\makebox(0,0)[bl]{\(e^{3,4}_{4}\) }}

\put(28,93){\oval(16,5)}


\put(49,91){\makebox(0,0)[bl]{\((1,3,0)\) }}

\put(00,93){\line(0,1){02}}

\put(04,93){\line(0,1){06}} \put(08,93){\line(0,1){06}}

\put(00,95){\line(1,0){8}} \put(04,97){\line(1,0){4}}
\put(04,99){\line(1,0){4}}

\put(00,93){\line(1,0){12}}

\put(00,91.5){\line(0,1){3}} \put(04,91.5){\line(0,1){3}}
\put(08,91.5){\line(0,1){3}} \put(12,91.5){\line(0,1){3}}

\put(01.5,89){\makebox(0,0)[bl]{\(1\)}}
\put(05.5,89){\makebox(0,0)[bl]{\(2\)}}
\put(09.5,89){\makebox(0,0)[bl]{\(3\)}}


\put(28,76){\line(0,1){14}}

\put(25,70.7){\makebox(0,0)[bl]{\(e^{3,4}_{5}\) }}

\put(28,73){\oval(16,5)}


\put(49,71){\makebox(0,0)[bl]{\((0,4,0)\) }}

\put(04,71.5){\line(0,1){08}} \put(08,71.5){\line(0,1){08}}

\put(04,73.5){\line(1,0){4}} \put(04,75.5){\line(1,0){4}}
\put(04,77.5){\line(1,0){4}} \put(04,79.5){\line(1,0){4}}

\put(00,71.5){\line(1,0){12}}

\put(00,70){\line(0,1){3}} \put(04,70){\line(0,1){3}}
\put(08,70){\line(0,1){3}} \put(12,70){\line(0,1){3}}

\put(01.5,67.5){\makebox(0,0)[bl]{\(1\)}}
\put(05.5,67.5){\makebox(0,0)[bl]{\(2\)}}
\put(09.5,67.5){\makebox(0,0)[bl]{\(3\)}}


\put(28,56){\line(0,1){14}}

\put(25,50.7){\makebox(0,0)[bl]{\(e^{3,4}_{6}\) }}

\put(28,53){\oval(16,5)}


\put(49,51){\makebox(0,0)[bl]{\((0,3,1)\) }}

\put(04,51){\line(0,1){06}} \put(08,51){\line(0,1){06}}
 \put(12,51){\line(0,1){02}}

\put(04,53){\line(1,0){8}} \put(04,55){\line(1,0){4}}
\put(04,57){\line(1,0){4}}

\put(00,50){\line(1,0){12}}

\put(00,49.5){\line(0,1){3}} \put(04,49.5){\line(0,1){3}}
\put(08,49.5){\line(0,1){3}} \put(12,49.5){\line(0,1){3}}

\put(01.5,47){\makebox(0,0)[bl]{\(1\)}}
\put(05.5,47){\makebox(0,0)[bl]{\(2\)}}
\put(09.5,47){\makebox(0,0)[bl]{\(3\)}}


\put(28,36){\line(0,1){14}}

\put(25,30.7){\makebox(0,0)[bl]{\(e^{3,4}_{7}\) }}

\put(28,33){\oval(16,5)}


\put(42,31){\makebox(0,0)[bl]{\((0,2,2)\) }}

\put(04,32){\line(0,1){04}}

\put(08,32){\line(0,1){04}} \put(12,32){\line(0,1){04}}
\put(04,34){\line(1,0){8}} \put(04,36){\line(1,0){8}}

\put(00,32){\line(1,0){12}}

\put(00,30.5){\line(0,1){3}} \put(04,30.5){\line(0,1){3}}
\put(08,30.5){\line(0,1){3}} \put(12,30.5){\line(0,1){3}}

\put(01.5,28){\makebox(0,0)[bl]{\(1\)}}
\put(05.5,28){\makebox(0,0)[bl]{\(2\)}}
\put(09.5,28){\makebox(0,0)[bl]{\(3\)}}


\put(28,24){\line(0,1){6}}

\put(25,18.7){\makebox(0,0)[bl]{\(e^{3,4}_{8}\) }}

\put(28,21){\oval(16,5)}


\put(42,19){\makebox(0,0)[bl]{\((0,1,3)\) }}

\put(04,21){\line(0,1){02}}

\put(08,21){\line(0,1){06}} \put(12,21){\line(0,1){06}}

\put(04,23){\line(1,0){8}} \put(08,25){\line(1,0){4}}
\put(08,27){\line(1,0){4}}

\put(00,21){\line(1,0){12}}

\put(00,19.5){\line(0,1){3}} \put(04,19.5){\line(0,1){3}}
\put(08,19.5){\line(0,1){3}} \put(12,19.5){\line(0,1){3}}

\put(01.5,17){\makebox(0,0)[bl]{\(1\)}}
\put(05.5,17){\makebox(0,0)[bl]{\(2\)}}
\put(09.5,17){\makebox(0,0)[bl]{\(3\)}}


\put(28,12){\line(0,1){6}}

\put(25,06.7){\makebox(0,0)[bl]{\(e^{3,4}_{12}\) }}

\put(28,09){\oval(16,5)}


\put(42,6.5){\makebox(0,0)[bl]{\((0,0,4)\) }}

\put(08,8.5){\line(0,1){08}} \put(12,8.5){\line(0,1){08}}
\put(08,10.5){\line(1,0){4}} \put(08,12.5){\line(1,0){4}}
\put(08,14.5){\line(1,0){4}} \put(08,16.5){\line(1,0){4}}

\put(00,8.5){\line(1,0){12}}

\put(00,07){\line(0,1){3}} \put(04,07){\line(0,1){3}}
\put(08,07){\line(0,1){3}} \put(12,07){\line(0,1){3}}

\put(01.5,4.5){\makebox(0,0)[bl]{\(1\)}}
\put(05.5,4.5){\makebox(0,0)[bl]{\(2\)}}
\put(09.5,4.5){\makebox(0,0)[bl]{\(3\)}}


\put(46.5,86.5){\line(-2,3){10.5}}

\put(45,80.7){\makebox(0,0)[bl]{\(e^{3,4}_{9}\) }}

\put(48,83){\oval(16,5)}


\put(57.5,81){\makebox(0,0)[bl]{\((2,1,1)\) }}

\put(00,84){\line(0,1){02}} \put(04,84){\line(0,1){04}}
\put(08,84){\line(0,1){04}} \put(12,84){\line(0,1){02}}

\put(00,86){\line(1,0){12}} \put(04,88){\line(1,0){4}}

\put(00,84){\line(1,0){12}}

\put(00,82.5){\line(0,1){3}} \put(04,82.5){\line(0,1){3}}
\put(08,82.5){\line(0,1){3}} \put(12,82.5){\line(0,1){3}}

\put(01.5,80){\makebox(0,0)[bl]{\(1\)}}
\put(05.5,80){\makebox(0,0)[bl]{\(2\)}}
\put(09.5,80){\makebox(0,0)[bl]{\(3\)}}


\put(44,66.5){\line(-1,2){11.4}}

\put(44,60){\line(-2,-1){8}}

\put(48,66){\line(0,1){14}}

\put(45,60.7){\makebox(0,0)[bl]{\(e^{3,4}_{10}\) }}

\put(48,63){\oval(16,5)}


\put(57.5,61){\makebox(0,0)[bl]{\((1,2,1)\) }}

\put(00,62){\line(0,1){02}} \put(04,62){\line(0,1){04}}
\put(08,62){\line(0,1){04}} \put(12,62){\line(0,1){02}}

\put(00,64){\line(1,0){12}} \put(04,66){\line(1,0){4}}

\put(00,62){\line(1,0){12}}

\put(00,60.5){\line(0,1){3}} \put(04,60.5){\line(0,1){3}}
\put(08,60.5){\line(0,1){3}} \put(12,60.5){\line(0,1){3}}

\put(01.5,58){\makebox(0,0)[bl]{\(1\)}}
\put(05.5,58){\makebox(0,0)[bl]{\(2\)}}
\put(09.5,58){\makebox(0,0)[bl]{\(3\)}}


\put(46.5,39.5){\line(-4,-1){12}}

\put(48,46){\line(0,1){14}}

\put(45,40.7){\makebox(0,0)[bl]{\(e^{3,4}_{11}\) }}

\put(48,43){\oval(16,5)}


\put(57.5,41){\makebox(0,0)[bl]{\((1,1,2)\) }}

\put(04,41.5){\line(0,1){04}}

\put(08,41.5){\line(0,1){04}} \put(12,41.5){\line(0,1){04}}
\put(04,43.5){\line(1,0){8}} \put(04,45.5){\line(1,0){8}}

\put(00,41.5){\line(1,0){12}}

\put(00,40){\line(0,1){3}} \put(04,40){\line(0,1){3}}
\put(08,40){\line(0,1){3}} \put(12,40){\line(0,1){3}}

\put(01.5,37.5){\makebox(0,0)[bl]{\(1\)}}
\put(05.5,37.5){\makebox(0,0)[bl]{\(2\)}}
\put(09.5,37.5){\makebox(0,0)[bl]{\(3\)}}


\end{picture}
\end{center}

 Further,  the problem is described as follows:
 \[ \max~ e(S),
  ~~~\max~ w(S),~~~
 s.t.
  ~~w(S) \geq 1  .\]
 As a result,
 we search for composite solutions
 which are nondominated by \(N(S)\)
 (i.e., Pareto-efficient).
 ``Maximization''  of \(e(S)\) is based on the corresponding poset.
 The considered combinatorial problem is NP-hard
 and an enumerative solving scheme is used.
%

 Further, combinatorial synthesis can be based on usage of
 interval multiset estimates of design alternatives for system parts.
%
 For the resultant system \(S = S(1) \star ... \star S(i) \star ... \star S(m) \)
 the same type of the interval multiset estimate is examined:
  an aggregated estimate (``generalized median'')
  of corresponding interval multiset estimates of its components
 (i.e., selected DAs).
 Thus, \( N(S) = (w(S);e(S))\), where
 \(e(S)\) is the ``generalized median'' of estimates of the solution
 components.
 Finally, the modified problem is:
%
 \[ \max~ e(S) = M^{g} =
  \arg \min_{M \in D}~~
  \sum_{i=1}^{m} ~ |\delta (M, e(S_{i})) |,
 ~~ \max~ w(S),~~~
 s.t.
%
  ~~w(S) \geq 1  .\]
 Here enumeration methods or heuristics are used
 (\cite{lev06},\cite{lev12morph},\cite{lev12a}).

 The next version of combinatorial synthesis problem will be the following.
 Let \(w(i,j)\) will be an interval multiset estimate
 for compatibility between two  selected design alternatives
 \(S_{i}\) and \(S_{j}\) in  solution \(S\)
 (the same type of interval multiset estimates is used for the elements
 and for their compatibility).
 Then, the problem is:
 \[ \max~ e(S) = M^{g} =
  \arg \min_{M \in D}~~
  \sum_{i=1}^{m} ~ |\delta (M, e(S_{i})) | + \sum_{(i,j) \in S} ~ |\delta (M, w(i,j))
  |,\]
%
%
 \[s.t.
%
  ~~~ w(S) = w(i,j) \succeq w_{0}, ~~ \forall (i,j) \in S  .\]
 Here \(w_{0}\) is a constraint
 (a bottom bound)
 for compatibility.
 Evidently,
 enumeration methods or heuristics can be  used here as well.

\subsection{Example of Composite Strategy}

 The four-stage (four-part)) morphological scheme of solving strategy
 for multicriteria ranking (sorting) to obtain the layered
 structure \(B\)
 is presented in
 Fig. 24:~
 \(S = H \star T \star U \star X\).
 Generally, design alternatives (local techniques, DAs)
 are evaluated upon six criteria  \cite{lev98}
 and the vector estimates are transformed into interval
 multiset estimates.
 In the example,
  expert judgment was used and illustrative estimates are presented
 in Table 4 (Fig. 27).
 Note, the estimates of DAs correspond to a certain applied
 situation,
 e.g., professional level of expert(s).

 Here the design of series solving strategy is based on the compressed
 morphology
 (without aggregation, i.e., without local techniques
 \(X_{1}\) and   \(X_{2}\))
 (Fig. 27).
 Table 5 contains estimates of compatibility.

 Finally, the following composite Pareto-efficient  DAs are obtained:

  (i)
   \(S_{1}=H_{1} \star T_{0} \star U_{2}\star X_{0}\), \(N(S_{1})=(2;4,0,0)\);
  (ii)
   \(S_{2}=H_{2} \star T_{0} \star U_{2}\star X_{0}\), \(N(S_{2})=(3;3,1,0)\);

  (iii)
   \(S_{3}=H_{3} \star T_{0} \star U_{2}\star X_{0}\), \(N(S_{3})=(3;3,1,0)\);
 (iv)
   \(S_{4}=H_{0} \star T_{0} \star U_{5}\star X_{0}\),
   \(N(S_{4})=(3;3,1,0)\).

 Fig. 28 illustrates a space of quality
 for the obtained solutions.

\begin{center}
\begin{picture}(83,94)

\put(17,92){\makebox(0,0)[bl]{Table 4. DAs and estimates}}

\put(00,00){\line(1,0){77}} \put(00,79){\line(1,0){77}}
\put(00,90){\line(1,0){77}}

\put(00,00){\line(0,1){90}} \put(08,00){\line(0,1){90}}
\put(62,00){\line(0,1){90}}

\put(77,00){\line(0,1){90}}

\put(63,86.5){\makebox(0,0)[bl]{Interval}}
\put(63,83.5){\makebox(0,0)[bl]{multiset}}
\put(63,80.5){\makebox(0,0)[bl]{estimate}}

\put(00.5,86){\makebox(0,0)[bl]{DAs}}
\put(09,85.6){\makebox(0,0)[bl]{Description}}


\put(01,74){\makebox(0,0)[bl]{\(H_{0}\)}}
\put(09,74.5){\makebox(0,0)[bl]{Absence}}

\put(68,73.5){\makebox(0,0)[bl]{\(-\)}}

\put(01,70){\makebox(0,0)[bl]{\(H_{1}\)}}
\put(09,70){\makebox(0,0)[bl]{Pairwise comparison}}

\put(64,69.5){\makebox(0,0)[bl]{\((4,0,0)\)}}

\put(01,66){\makebox(0,0)[bl]{\(H_{2}\)}}
\put(09,66){\makebox(0,0)[bl]{Dominance by Pareto-rule}}

\put(64,65.5){\makebox(0,0)[bl]{\((3,1,0)\)}}

\put(01,62){\makebox(0,0)[bl]{\(H_{3}\)}}
\put(09,62){\makebox(0,0)[bl]{Outranking technique}}

\put(64,61.5){\makebox(0,0)[bl]{\((3,1,0)\)}}


\put(01,58){\makebox(0,0)[bl]{\(T_{0}\)}}
\put(09,58.5){\makebox(0,0)[bl]{Absence}}

\put(68,57.5){\makebox(0,0)[bl]{\(-\)}}

\put(01,54){\makebox(0,0)[bl]{\(T_{1}\)}}
\put(09,54){\makebox(0,0)[bl]{Line elements sum of preference}}
\put(09,50.7){\makebox(0,0)[bl]{matrix}}

\put(64,53.5){\makebox(0,0)[bl]{\((1,2,1)\)}}

\put(01,46){\makebox(0,0)[bl]{\(T_{2}\)}}
\put(09,46){\makebox(0,0)[bl]{Additive utility function}}

\put(64,45.5){\makebox(0,0)[bl]{\((2,2,0)\)}}


\put(01,42){\makebox(0,0)[bl]{\(U_{0}\)}}
\put(09,42.5){\makebox(0,0)[bl] {Absence}}

\put(68,41.5){\makebox(0,0)[bl]{\(-\)}}

\put(01,38){\makebox(0,0)[bl]{\(U_{1}\)}}
\put(09,38){\makebox(0,0)[bl] {Step-by-step revelation }}

\put(09,34.5){\makebox(0,0)[bl] {{\it of maximal} elements}}

\put(64,37.5){\makebox(0,0)[bl]{\((2,2,0)\)}}

\put(01,30){\makebox(0,0)[bl]{\(U_{2}\)}}
\put(09,30){\makebox(0,0)[bl]{Step-by-step revelation of Pareto}}

\put(09,26.5){\makebox(0,0)[bl]{efficient elements}}

\put(64,29.5){\makebox(0,0)[bl]{\((4,0,0)\)}}

\put(01,22){\makebox(0,0)[bl]{\(U_{3}\)}}
\put(09,22){\makebox(0,0)[bl]{Dividing the linear ranking}}

\put(64,21.5){\makebox(0,0)[bl]{\((0,1,3)\)}}

\put(01,18){\makebox(0,0)[bl]{\(U_{4}\)}}
\put(09,18){\makebox(0,0)[bl]{Expert procedure for ranking}}

\put(64,17.5){\makebox(0,0)[bl]{\((2,2,0)\)}}

\put(01,14){\makebox(0,0)[bl]{\(U_{5}\)}}
\put(09,14){\makebox(0,0)[bl]{Direct logical method for ranking}}

\put(64,13.5){\makebox(0,0)[bl]{\((3,1,0)\)}}


\put(01,10){\makebox(0,0)[bl]{\(X_{0}\)}}
\put(09,10.5){\makebox(0,0)[bl]{Absence}}

\put(64,09.5){\makebox(0,0)[bl]{\((0,4,0)\)}}

\put(01,06){\makebox(0,0)[bl]{\(X_{1}\)}}
\put(09,06){\makebox(0,0)[bl]{Expert procedure for ranking}}

\put(64,05.5){\makebox(0,0)[bl]{\((3,1,0)\)}}

\put(01,02){\makebox(0,0)[bl]{\(X_{2}\)}}
\put(09,02){\makebox(0,0)[bl]{Direct logical method for ranking}}
\put(64,1.5){\makebox(0,0)[bl]{\((4,0,0)\)}}

\end{picture}
%
\begin{picture}(62,60)

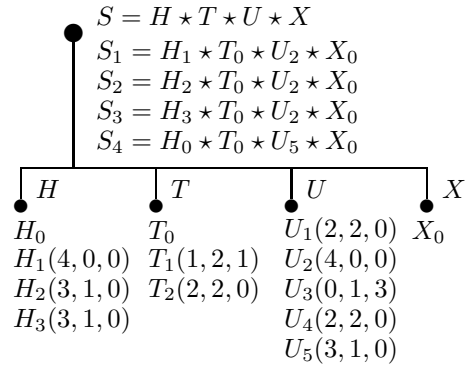
\put(00,00){\makebox(0,0)[bl]{Fig. 27. Morphology of series
 strategy}}


\put(1,28){\line(0,1){4}} \put(01,27){\circle*{1.6}}
\put(03,28){\makebox(0,0)[bl]{\(H\)}}

\put(00,22){\makebox(0,0)[bl]{\(H_{0}\)}}
\put(00,18){\makebox(0,0)[bl]{\(H_{1}(4,0,0)\)}}
\put(00,14){\makebox(0,0)[bl]{\(H_{2}(3,1,0)\)}}
\put(00,10){\makebox(0,0)[bl]{\(H_{3}(3,1,0)\)}}


\put(19,28){\line(0,1){4}} \put(19,27){\circle*{1.6}}
\put(21,28){\makebox(0,0)[bl]{\(T\)}}

\put(18,22){\makebox(0,0)[bl]{\(T_{0}\)}}
\put(18,18){\makebox(0,0)[bl]{\(T_{1}(1,2,1)\)}}
\put(18,14){\makebox(0,0)[bl]{\(T_{2}(2,2,0)\)}}


\put(37,28){\line(0,1){4}} \put(37,27){\circle*{1.6}}
\put(39,28){\makebox(0,0)[bl]{\(U\)}}

\put(36,22){\makebox(0,0)[bl]{\(U_{1}(2,2,0)\)}}
\put(36,18){\makebox(0,0)[bl]{\(U_{2}(4,0,0)\)}}
\put(36,14){\makebox(0,0)[bl]{\(U_{3}(0,1,3)\)}}
\put(36,10){\makebox(0,0)[bl]{\(U_{4}(2,2,0)\)}}
\put(36,06){\makebox(0,0)[bl]{\(U_{5}(3,1,0)\)}}


\put(55,28){\line(0,1){4}} \put(55,27){\circle*{1.6}}
\put(57,28){\makebox(0,0)[bl]{\(X\)}}

\put(53,22){\makebox(0,0)[bl]{\(X_{0}\)}}


\put(01,32){\line(1,0){54}} \put(08,32){\line(0,1){18}}
\put(08,50){\circle*{2.5}}

\put(11,51){\makebox(0,0)[bl]{\(S = H \star T \star U \star X \)}}

\put(11,46){\makebox(0,0)[bl] {\(S_{1}= H_{1} \star T_{0} \star
U_{2} \star X_{0} \)}}

\put(11,42){\makebox(0,0)[bl] {\(S_{2}= H_{2} \star T_{0} \star
U_{2} \star X_{0} \)}}

\put(11,38){\makebox(0,0)[bl] {\(S_{3}= H_{3} \star T_{0} \star
U_{2} \star X_{0} \)}}

\put(11,34){\makebox(0,0)[bl] {\(S_{4}= H_{0} \star T_{0} \star
U_{5} \star X_{0} \)}}

\end{picture}
\end{center}

\begin{center}
\begin{picture}(62,63)

\put(02.5,58){\makebox(0,0)[bl]{Table 5. Compatibility of DAs}}

\put(00,56){\line(1,0){52}} \put(00,50){\line(1,0){52}}
\put(00,00){\line(1,0){52}}

\put(00,00){\line(0,1){56}} \put(06,00){\line(0,1){56}}
\put(52,00){\line(0,1){56}}

\put(11,50){\line(0,1){06}} \put(16,50){\line(0,1){06}}
\put(21,50){\line(0,1){06}} \put(26,50){\line(0,1){06}}
\put(31,50){\line(0,1){06}} \put(36,50){\line(0,1){06}}
\put(41,50){\line(0,1){06}} \put(46,50){\line(0,1){06}}

\put(07,52){\makebox(0,8)[bl]{\(T_{0}\)}}
\put(12,52){\makebox(0,8)[bl]{\(T_{1}\)}}
\put(17,52){\makebox(0,8)[bl]{\(T_{2}\)}}

\put(22,52){\makebox(0,8)[bl]{\(U_{1}\)}}
\put(27,52){\makebox(0,8)[bl]{\(U_{2}\)}}
\put(32,52){\makebox(0,8)[bl]{\(U_{3}\)}}
\put(37,52){\makebox(0,8)[bl]{\(U_{4}\)}}
\put(42,52){\makebox(0,8)[bl]{\(U_{5}\)}}

\put(47,52){\makebox(0,8)[bl]{\(X_{0}\)}}


\put(01,46){\makebox(0,8)[bl]{\(H_{0}\)}}
\put(08,46){\makebox(0,8)[bl]{\(3\)}}
\put(13,46){\makebox(0,8)[bl]{\(3\)}}
\put(18,46){\makebox(0,8)[bl]{\(3\)}}

\put(23,46){\makebox(0,8)[bl]{\(0\)}}
\put(28,46){\makebox(0,8)[bl]{\(0\)}}
\put(33,46){\makebox(0,8)[bl]{\(0\)}}
\put(38,46){\makebox(0,8)[bl]{\(3\)}}
\put(43,46){\makebox(0,8)[bl]{\(3\)}}
\put(48,46){\makebox(0,8)[bl]{\(0\)}}

\put(01,42){\makebox(0,8)[bl]{\(H_{1}\)}}
\put(08,42){\makebox(0,8)[bl]{\(0\)}}
\put(13,42){\makebox(0,8)[bl]{\(1\)}}
\put(18,42){\makebox(0,8)[bl]{\(0\)}}

\put(23,42){\makebox(0,8)[bl]{\(2\)}}
\put(28,42){\makebox(0,8)[bl]{\(2\)}}
\put(33,42){\makebox(0,8)[bl]{\(0\)}}
\put(38,42){\makebox(0,8)[bl]{\(0 \)}}
\put(43,42){\makebox(0,8)[bl]{\(0\)}}
\put(48,42){\makebox(0,8)[bl]{\(3\)}}

\put(01,38){\makebox(0,8)[bl]{\(H_{2}\)}}
\put(08,38){\makebox(0,8)[bl]{\(0\)}}
\put(13,38){\makebox(0,8)[bl]{\(1\)}}
\put(18,38){\makebox(0,8)[bl]{\(0\)}}

\put(23,38){\makebox(0,8)[bl]{\(3\)}}
\put(28,38){\makebox(0,8)[bl]{\(3\)}}
\put(33,38){\makebox(0,8)[bl]{\(0\)}}
\put(38,38){\makebox(0,8)[bl]{\(0\)}}
\put(43,38){\makebox(0,8)[bl]{\(0\)}}
\put(48,38){\makebox(0,8)[bl]{\(3\)}}

\put(01,34){\makebox(0,8)[bl]{\(H_{3}\)}}
\put(08,34){\makebox(0,8)[bl]{\(0\)}}
\put(13,34){\makebox(0,8)[bl]{\(1\)}}
\put(18,34){\makebox(0,8)[bl]{\(0\)}}

\put(23,34){\makebox(0,8)[bl]{\(3\)}}
\put(28,34){\makebox(0,8)[bl]{\(3\)}}
\put(33,34){\makebox(0,8)[bl]{\(0\)}}
\put(38,34){\makebox(0,8)[bl]{\(0\)}}
\put(43,34){\makebox(0,8)[bl]{\(0\)}}
\put(48,34){\makebox(0,8)[bl]{\(3\)}}


\put(01,30){\makebox(0,8)[bl]{\(T_{0}\)}}
\put(23,30){\makebox(0,8)[bl]{\(0\)}}
\put(28,30){\makebox(0,8)[bl]{\(0\)}}
\put(33,30){\makebox(0,8)[bl]{\(0\)}}
\put(38,30){\makebox(0,8)[bl]{\(3\)}}
\put(43,30){\makebox(0,8)[bl]{\(3\)}}
\put(48,30){\makebox(0,8)[bl]{\(0\)}}

\put(01,26){\makebox(0,8)[bl]{\(T_{1}\)}}
\put(23,26){\makebox(0,8)[bl]{\(0\)}}
\put(28,26){\makebox(0,8)[bl]{\(0\)}}
\put(33,26){\makebox(0,8)[bl]{\(2\)}}
\put(38,26){\makebox(0,8)[bl]{\(0\)}}
\put(43,26){\makebox(0,8)[bl]{\(0\)}}
\put(48,26){\makebox(0,8)[bl]{\(3\)}}

\put(01,22){\makebox(0,8)[bl]{\(T_{2}\)}}
\put(23,22){\makebox(0,8)[bl]{\(0\)}}
\put(28,22){\makebox(0,8)[bl]{\(0\)}}
\put(33,22){\makebox(0,8)[bl]{\(2\)}}
\put(38,22){\makebox(0,8)[bl]{\(0\)}}
\put(43,22){\makebox(0,8)[bl]{\(0\)}}
\put(48,22){\makebox(0,8)[bl]{\(3\)}}


\put(01,18){\makebox(0,8)[bl]{\(U_{1}\)}}
\put(48,18){\makebox(0,8)[bl]{\(3\)}}

\put(01,14){\makebox(0,8)[bl]{\(U_{2}\)}}
\put(48,14){\makebox(0,8)[bl]{\(3\)}}

\put(01,10){\makebox(0,8)[bl]{\(U_{3}\)}}
\put(48,10){\makebox(0,8)[bl]{\(3\)}}

\put(01,06){\makebox(0,8)[bl]{\(U_{4}\)}}
\put(48,06){\makebox(0,8)[bl]{\(3\)}}

\put(01,02){\makebox(0,8)[bl]{\(U_{5}\)}}
\put(48,02){\makebox(0,8)[bl]{\(3\)}}

\end{picture}
%
\begin{picture}(70,66)
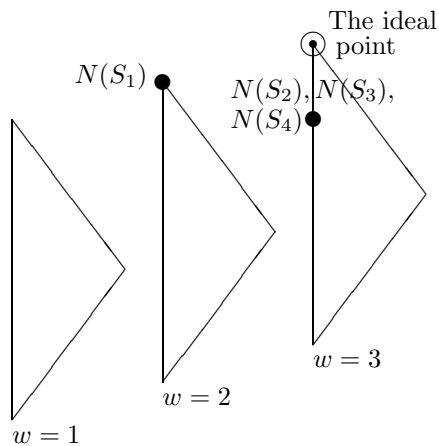
\put(11.5,00){\makebox(0,0)[bl]{Fig. 28. Space of system quality}}

\put(10,010){\line(0,1){40}} \put(10,010){\line(3,4){15}}
\put(10,050){\line(3,-4){15}}

\put(30,015){\line(0,1){40}} \put(30,015){\line(3,4){15}}
\put(30,055){\line(3,-4){15}}

\put(50,020){\line(0,1){40}} \put(50,020){\line(3,4){15}}
\put(50,060){\line(3,-4){15}}

\put(30,55){\circle*{2}}
\put(18.3,54){\makebox(0,0)[bl]{\(N(S_{1})\)}}

\put(50,50){\circle*{2}}
\put(39,52){\makebox(0,0)[bl]{\(N(S_{2}),N(S_{3}),\)}}
\put(39,48){\makebox(0,0)[bl]{\(N(S_{4})\)}}

\put(50,60){\circle*{1}} \put(50,60){\circle{3}}

\put(52,62){\makebox(0,0)[bl]{The ideal}}
\put(53,58){\makebox(0,0)[bl]{point}}

\put(10,7){\makebox(0,0)[bl]{\(w=1\)}}
\put(30,12){\makebox(0,0)[bl]{\(w=2\)}}
\put(50,17){\makebox(0,0)[bl]{\(w=3\)}}
\end{picture}
\end{center}

\section{CONCLUSION}

 In the paper,
 a modular approach to
 solving strategies in DSS for multicriteria ranking (sorting)
 has been described.
 First, methodological issues
 in architectural design of DSS are examined
 (decision making framework, typical problems,
 approaches to configuration of solving strategies).
 Second, combinatorial synthesis
 for design of series composite strategy
 for multicriteria ranking (sorting))
 is presented.
 The study is based on DSS COMBI
  for multicriteria ranking (1984...1991).
 Evidently,
 an important and prospective direction consists in
 aggregation (fusion) of preference relations which are obtained
 from different sources (experts, algorithm / procedures).
 This leads to series-parallel solving strategies
 (e.g., \cite{lev98}).
 Here it may be reasonable
 to take into account the following two research directions:
 (1) routing in And-Or graphs
 (e.g., \cite{ade02},\cite{din90});
 (2) activity nets \cite{elm95}).
 In addition,
 it is reasonable to point out the significance
 of an extended composite solving framework
 that involves
 a preliminary stage for problem formulation,
 a final stage for analysis of results,
 and
 on-line monitoring of the solving process
 (Fig. 28).
 Here, the same design approaches may be applied.

\begin{center}
\begin{picture}(137,91)
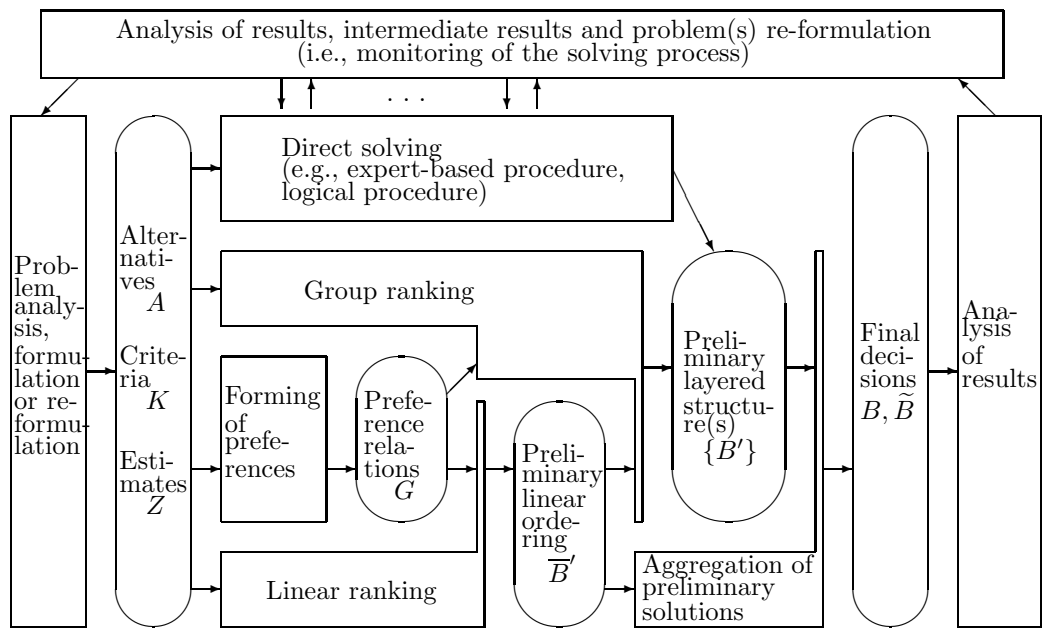
\put(40.5,00){\makebox(0,0)[bl] {Fig. 28. Extended functional
  graph}}


\put(4,80){\line(0,1){09}} \put(132,80){\line(0,1){09}}
\put(4,80){\line(1,0){128}} \put(4,89){\line(1,0){128}}

\put(14,84.5){\makebox(0,0)[bl]{Analysis of results, intermediate
results and problem(s) re-formulation}}

\put(38,81.5){\makebox(0,0)[bl]{(i.e., monitoring of the solving
process)}}

\put(50,77){\makebox(0,0)[bl]{.~.~.}}

\put(36,80){\vector(0,-1){4}} \put(40,76){\vector(0,1){4}}
\put(66,80){\vector(0,-1){4}} \put(70,76){\vector(0,1){4}}

\put(09,80){\vector(-1,-1){5}} \put(131,75){\vector(-1,1){5}}


\put(00,07){\line(0,1){68}} \put(10,07){\line(0,1){68}}
\put(00,07){\line(1,0){10}} \put(00,75){\line(1,0){10}}

\put(00.5,54){\makebox(0,0)[bl]{Prob-}}
\put(00.5,51){\makebox(0,0)[bl]{lem}}

\put(00.5,48){\makebox(0,0)[bl]{analy-}}
\put(00.5,45){\makebox(0,0)[bl]{sis,}}

\put(00.5,42){\makebox(0,0)[bl]{formu-}}
\put(00.5,39){\makebox(0,0)[bl]{lation}}
\put(00.5,36){\makebox(0,0)[bl]{or re-}}
\put(00.5,33){\makebox(0,0)[bl]{formu-}}
\put(00.5,30){\makebox(0,0)[bl]{lation}}

\put(10,41){\vector(1,0){4}}


\put(19,41){\oval(10,68)}

\put(14.5,58){\makebox(0,0)[bl]{Alter-}}
\put(14.5,55){\makebox(0,0)[bl]{nati-}}
\put(14.5,52){\makebox(0,0)[bl]{ves}}
\put(18,49){\makebox(0,0)[bl]{\(A\)}}

\put(14.5,42){\makebox(0,0)[bl]{Crite-}}
\put(14.5,39){\makebox(0,0)[bl]{ria}}
\put(18,36){\makebox(0,0)[bl]{\(K\)}}

\put(14.5,28){\makebox(0,0)[bl]{Esti-}}
\put(14.5,25){\makebox(0,0)[bl]{mates}}
\put(18,22){\makebox(0,0)[bl]{\(Z\)}}


\put(117,41){\oval(10,68)}

\put(113,45){\makebox(0,0)[bl]{Final}}
\put(113,42){\makebox(0,0)[bl]{deci-}}
\put(113,39){\makebox(0,0)[bl]{sions}}
\put(113,34){\makebox(0,0)[bl]{\(B,\widetilde{B}\)}}

\put(108,28){\vector(1,0){4}}


\put(126,07){\line(0,1){68}} \put(137,07){\line(0,1){68}}
\put(126,07){\line(1,0){11}} \put(126,75){\line(1,0){11}}

\put(126.5,48){\makebox(0,0)[bl]{Ana-}}
\put(126.5,45){\makebox(0,0)[bl]{lysis}}
\put(126.5,42){\makebox(0,0)[bl]{of}}
\put(126.5,39){\makebox(0,0)[bl]{results}}

\put(122,41){\vector(1,0){4}}


\put(62,37){\line(1,0){1}}

\put(62,17){\line(0,1){20}}

\put(28,07){\line(0,1){10}} \put(63,07){\line(0,1){30}}

\put(28,17){\line(1,0){34}} \put(28,07){\line(1,0){35}}

\put(34,10){\makebox(0,0)[bl]{Linear ranking}}

\put(24,12){\vector(1,0){4}}


\put(28,21){\line(0,1){22}} \put(42,21){\line(0,1){22}}
\put(28,21){\line(1,0){14}} \put(28,43){\line(1,0){14}}

\put(28.5,36){\makebox(0,0)[bl]{Forming}}
\put(28.5,33){\makebox(0,0)[bl]{of}}
\put(28.5,30){\makebox(0,0)[bl]{prefe-}}
\put(28.5,27){\makebox(0,0)[bl]{rences}}

\put(24,28){\vector(1,0){4}}

\put(42,28){\vector(1,0){4}}

\put(58,28){\vector(1,0){4}}

\put(63,28){\vector(1,0){4}}


\put(52,32){\oval(12,22)}

\put(47,36){\makebox(0,0)[bl]{Prefe-}}
\put(47,33){\makebox(0,0)[bl]{rence}}
\put(47,30){\makebox(0,0)[bl]{rela-}}
\put(47,27){\makebox(0,0)[bl]{tions}}
\put(51,24){\makebox(0,0)[bl]{\(G\)}}


\put(73,22){\oval(12,30)}

\put(68,29){\makebox(0,0)[bl]{Preli-}}
\put(68,25.6){\makebox(0,0)[bl]{minary}}
\put(68,23){\makebox(0,0)[bl]{linear}}
\put(68,20){\makebox(0,0)[bl]{orde-}}
\put(68,17){\makebox(0,0)[bl]{ring}}
\put(71.5,13){\makebox(0,0)[bl]{\(\overline{B}'\)}}

\put(79,28){\vector(1,0){4}}


\put(107,57){\line(1,0){1}}

\put(107,17){\line(0,1){40}}

\put(83,07){\line(0,1){10}} \put(108,07){\line(0,1){50}}

\put(83,17){\line(1,0){24}} \put(83,07){\line(1,0){25}}

\put(84,13.5){\makebox(0,0)[bl]{Aggregation of}}
\put(84,10.6){\makebox(0,0)[bl]{preliminary}}
\put(84,08){\makebox(0,0)[bl]{solutions}}

\put(79,12){\vector(1,0){4}}



\put(28,47){\line(0,1){10}}

\put(83,40){\line(-1,0){21}}

\put(62,40){\line(0,1){7}}

\put(83,21){\line(0,1){19}} \put(84,21){\line(0,1){36}}
\put(83,21){\line(1,0){1}}

\put(28,47){\line(1,0){34}} \put(28,57){\line(1,0){56}}

\put(39,50){\makebox(0,0)[bl]{Group ranking}}

\put(24,52){\vector(1,0){4}}

\put(58,38){\vector(1,1){4}}

\put(84,41.5){\vector(1,0){4}}


\put(28,61){\line(0,1){14}} \put(88,61){\line(0,1){14}}

\put(28,61){\line(1,0){60}} \put(28,75){\line(1,0){60}}

\put(36,69){\makebox(0,0)[bl]{Direct solving}}
\put(36,66){\makebox(0,0)[bl]{(e.g., expert-based procedure,}}
\put(36,63){\makebox(0,0)[bl]{logical procedure)}}

\put(24,68){\vector(1,0){4}} \put(88,68){\vector(1,-2){5.5}}



\put(95.5,39){\oval(15,36)}

\put(89.5,44){\makebox(0,0)[bl]{Preli-}}
\put(89.5,41){\makebox(0,0)[bl]{minary}}
\put(89.5,38){\makebox(0,0)[bl]{layered}}
\put(89.5,35.4){\makebox(0,0)[bl]{structu-}}
\put(89.5,32){\makebox(0,0)[bl]{re(s)}}
\put(92,28.5){\makebox(0,0)[bl]{\(\{B'\}\)}}

\put(103,41.5){\vector(1,0){4}}

\end{picture}
\end{center}

 Finally, it is reasonable to emphasize the significance of our
    approach for teaching of multicriteria analysis and system
    design
    (techniques, case studies, applied examples, projects)
    (e.g., \cite{lev95e},\cite{leved00},\cite{lev06edu},\cite{lev10edu}).

\end{document}